\documentclass[aos,MSNbibl,seceqn,dvips]{arximspdf}
\usepackage{graphicx}

\doi{10.1214/12-AOS1045} 
\volume{40}
\issue{5}
\pubyear{2012}
\firstpage{2634}
\lastpage{2666}

\makeatletter

\newcommand{\hhat}[1]{\hspace*{1.5pt}\hat{\hspace*{-1.5pt}#1}}

\newproclaim{Remark}{Remark}
\newproclaim{Remarks}{Remarks}

\newproclaim{step}{Step}
\newproclaim{example}{Example}
\newproclaim{assumption}{Assumption}

\newtheorem{theorem}{Theorem}
\newtheorem{corollary}{Corollary}
\newtheorem{lemma}{Lemma}

\newcommand{\rrvert}{\vert}
\newcommand{\llvert}{\vert}

\newcommand{\dddot}[1]{\hspace*{2pt}\dot{\vphantom{B}\hspace*{2pt}}\ddot{\hspace*{-3pt}#1}}

\newcommand{\ve}[1]{\bolds{#1}}

\makeatother

\begin{document}
\begin{frontmatter}

\title{Multivariate varying coefficient model for functional responses}
\runtitle{Multivariate varying coefficient model}

\begin{aug}
\author[A]{\fnms{Hongtu} \snm{Zhu}\thanksref{t1}\ead[label=e1]{hzhu@bios.unc.edu}},
\author[B]{\fnms{Runze} \snm{Li}\corref{}\thanksref{t2}\ead[label=e2]{rli@stat.psu.edu}}
\and
\author[C]{\fnms{Linglong} \snm{Kong}\thanksref{t1}\ead[label=e3]{lkong@ualberta.ca}}
\runauthor{H. Zhu, R. Li and L. Kong}
\affiliation{University of North Carolina at Chapel Hill,
Pennsylvania State University and~University of Alberta}
\address[A]{H. Zhu\\
Departments of Biostatistics\\
Biomedical Research Imaging Center \\
University of North Carolina at Chapel Hill\\
Chapel Hill, North Carolina 27599\\
USA\\
\printead{e1}}
\address[B]{R. Li\\
Department of Statistics\\
Pennsylvania State University\\
University Park, Pennsylvania 16802\\
USA\\
\printead{e2}}
\address[C]{L. Kong\\
Department of Mathematical\\
\quad and Statistical Sciences\\
University of Alberta\\
CAB 632\\
Edmonton, Alberta, T6G 2G1\\
Canada\\
\printead{e3}} 
\end{aug}

\thankstext{t1}{Supported by NIH Grants RR025747-01, P01CA142538-01, MH086633,
EB005149-01 and AG033387.}

\thankstext{t2}{Supported by NSF Grant DMS-03-48869, NIH Grants
P50-DA10075 and
R21-DA024260 and NNSF of China 11028103. The content is solely the
responsibility of the authors and does not necessarily represent the
official views of the NSF or the NIH.}

\received{\smonth{3} \syear{2012}}
\revised{\smonth{9} \syear{2012}}

%
\begin{abstract}
Motivated by recent work studying massive imaging data in the
neuroimaging literature, we propose multivariate varying coefficient
models (MVCM) for modeling the relation between multiple functional
responses and a set of covariates. We develop several statistical
inference procedures for MVCM and systematically study their
theoretical properties. We first establish the weak convergence of
the local linear estimate of coefficient functions, as well as
its asymptotic bias and variance, and then we derive asymptotic bias
and mean integrated squared error of smoothed individual functions and
their uniform convergence rate. We establish the uniform convergence
rate of the estimated covariance function of the individual functions
and its associated eigenvalue and eigenfunctions. We propose a global
test for linear hypotheses of varying coefficient functions, and
derive its asymptotic distribution under the null hypothesis. We also
propose a simultaneous confidence band for each individual effect
curve. We conduct Monte Carlo simulation to examine the finite-sample
performance of the proposed procedures. We apply MVCM to investigate
the development of white matter diffusivities along the genu tract of
the corpus callosum in a clinical study of neurodevelopment.
\end{abstract}

%

\end{frontmatter}

\section{Introduction}\label{sec1}


With modern imaging techniques, massive imaging data can be
observed over both time and space \cite
{Towle1993,Fass2008,Niedermeyer2004,Buzsaki2006,Friston2009,Heywood2006}. Such
imaging techniques include functional magnetic resonance imaging
(fMRI), electroencephalography (EEG),\vadjust{\goodbreak} diffusion tensor imaging
(DTI), positron emission tomography (PET) and single photon
emission-computed tomography (SPECT) among many other imaging
techniques. See, for example, a recent review of multiple
biomedical imaging techniques and their applications in cancer
detection and prevention in Fass~\cite{Fass2008}. Among them,
predominant functional imaging techniques including fMRI and EEG
have been widely used in behavioral and cognitive neuroscience to
understand functional segregation and integration of different
brain regions in a single subject and across different
populations~\cite{Friston2009,Friston2007,Huettel2004}. In DTI,
multiple diffusion properties are measured along common major white
matter fiber tracts across multiple subjects to characterize the
structure and orientation of white matter structure in human brain
in vivo~\cite{Basser1994b,Basser1994a,Zhu2007b}.

\begin{figure}[b]
\begin{tabular}{@{}cc@{}}

\includegraphics{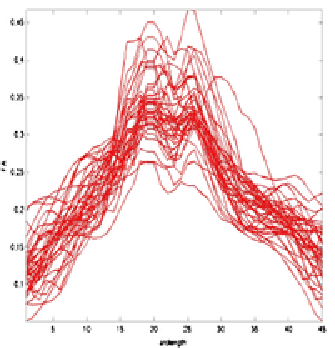}
 & \includegraphics{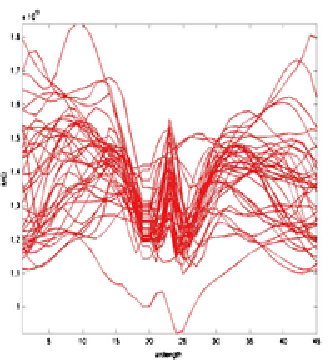}\\
(a) & (b)\\[6pt]

\includegraphics{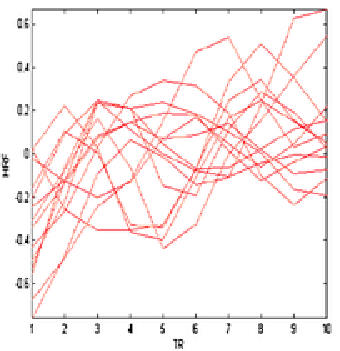}
 & \includegraphics{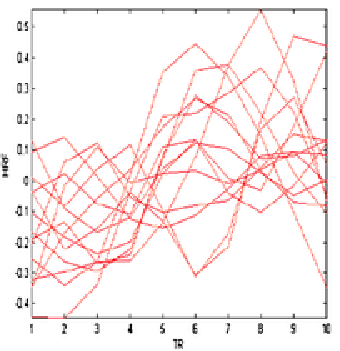}\\
(c) & (d)
\end{tabular}
\caption{Representative functional neuroimaging data: \textup{(a)} and \textup{(b)} FA
and MD along the genu tract of the corpus callosum from 40
randomly selected infants; and \textup{(c)} and~\textup{(d)} the estimated hemodynamic
response functions (HRF) corresponding to two stimulus categories from
14 subjects.} \label{tract3fig0}
\end{figure}

A common feature of many imaging techniques is that massive
functional data are observed/calculated at the same design points, such
as time for functional images (e.g., PET and fMRI).
As an illustration, we present two smoothed functional data as an
illustration and a real imaging data in Section~\ref{sec6}, that we encounter in
neuroimaging studies.
First, we plot two diffusion properties, called fractional
anisotropy (FA) and mean diffusivity (MD), measured at 45 grid points
along the genu tract of the corpus callosum
[Figure~\ref{tract3fig0}(a) and (b)] from 40 randomly selected
infants from a clinical study of neurodevelopment with more than 500
infants. Scientists are
particularly interested in delineating the structure of the
variability of these functional FA and MD data and their association
with a set of covariates of interest, such
as age. We will systematically investigate the development of FA and MD
along the genu of the corpus callosum tract in Section~\ref{sec6}.
Second,
we consider
the BOLD fMRI signal, which is based on hemodynamic responses secondary
to neural activity. We
plot the estimated hemodynamic response functions (HRF) corresponding
to two stimulus categories from 14 randomly selected subjects at a
selected voxel of a common template space from a clinical study of
Alzheimer's disease with more than 100 infants.
Although the canonical form of the HRF is often used, when applying
fMRI in a clinical
population with possibly altered hemodynamic responses [Figure~\ref{tract3fig0}(c)
and~(d)], using the subject's own HRF in
fMRI data analysis may be advantageous because HRF variability is
greater across subjects
than across brain regions within a subject \cite
{Lindquist2008,AguirreZarahn1998}. We
are particularly interested in delineating the structure of the
variability of the HRF and their association with a set of covariates
of interest, such
as diagnostic group~\cite{Lindquist2008b}.

A varying-coefficient
model, which allows its regression coefficients to vary over some
predictors of interest, is a powerful statistical tool for
addressing these scientific questions.
Since it was systematically introduced to statistical literature
by Hastie and Tibshirani~\cite{HastieTibshirani1993}, many
varying-coefficient models have
been widely studied and developed for longitudinal, time series
and functional data \cite
{MR1742497,MR1769751,MR1959093,MR2425354,MR2504204,MR1666699,Ramsay2005,MR2087972,MR1888349,Zhang2007,MR2523900}.
However, most varying-coefficient models in the existing literature
are developed for univariate response.
Let ${\mathbf y}_{i}(s)=(y_{i1}(s),\ldots, y_{i J}(s))^T$ be a
$J$-dimensional functional response vector for subject $i$, $i=1,\ldots,
n$, and ${\mathbf x}_i$ be its associated $p\times1$ vector of
covariates of interest. Moreover, $s$ varies in a compact subset of
Euclidean space and
denotes the design point, such as time for functional images and voxel
for structural and functional images. For notational simplicity, we
assume $s\in[0, 1]$, but our results
can be easily extended to higher dimensions.
A \textit{multivariate varying coefficient model} (\textit{MVCM}) is defined as
%
\begin{equation}
\label{Tract3eq1} {y}_{ij}(s)= {\mathbf x}_i^TB_j(s)+
\eta_{ij}(s)+\varepsilon_{ij}(s) \qquad\mbox{for }j=1,\ldots, J,
\end{equation}
where $B_j(s)=(b_{j1}(s),\ldots, b_{jp}(s))^T$ is a $p\times1$
vector of functions of $s$, $\varepsilon_{ij}(s)$ are measurement
errors, and $\eta_{ij}(s)$ characterizes individual curve
variations from ${\mathbf x}_i^TB_j(s)$. Moreover, $\{\eta_{i j}(s)\dvtx
s\in[0, 1]\}$ is assumed to be a stochastic process indexed by
$s\in[0, 1]$ and used to characterize the within-curve
dependence. For image data, it is typical that the $J$ functional
responses ${\mathbf y}_{i}(s)$ are measured at the same location for
all subjects and exhibit both the within-curve and between-curve
dependence structure. Thus, for ease of notation, it is assumed
throughout this paper that
${\mathbf y}_{i}(s)$ was measured at the same $M$ location points
$s_1=0\leq s_2\leq\cdots\leq s_{M}=1$ for all $i$.


Most varying coefficient models in the existing literature coincide
model (\ref{Tract3eq1}) with $J=1$ and without the within-curve
dependence.
Statistical inferences for these varying coefficient models have been
relatively well
studied. Particularly, Hoover et al.~\cite{MR1666699} and Wu, Chiang
and Hoover~\cite{WuChiang1998} were
among the first to introduce the time-varying coefficient models for
analysis of longitudinal data.
Recently, Fan and Zhang~\cite{MR2425354} gave a
comprehensive review of
various statistical procedures proposed for many varying coefficient
models. It is of particular
interest in data analysis to construct simultaneous confidence bands
(SCB) for
any linear combination of $B_j$ instead of pointwise confidence
intervals and to
develop global test statistics for the general hypothesis testing
problem on $B_j$.
For univariate varying coefficient models without the
within-curve dependence,
Fan and Zhang~\cite{MR1804172} constructed SCB using the limit theory
for the maximum of the normalized deviation of the estimate from its
expected value. Faraway~\cite{Faraway1997}, Chiou, M{\"u}ller and
Wang~\cite {Chiou2004}, and Cardot~\cite{Cardot2007} proposed several
varying coefficient models and their associated estimators for
univariate functional response, but they did not give functional
central limit theorem and simultaneous confidence band for their
estimators. It has been technically difficult to carry out statistical
inferences including simultaneous confidence band and global test
statistic on $B_j$ in the presence of the within-curve dependence.


There have been several recent attempts to solve this problem in
various settings.
For time series data, which may be viewed as a case with $n=1$ and
$M\rightarrow\infty$, asymptotic SCB
for coefficient functions in varying coefficient models can be built
by using local kernel regression and a Gaussian approximation
result for nonstationary time series~\cite{ZhouWu2010}. For sparse
irregular longitudinal data,
Ma, Yang and Carroll~\cite{MaYang2011} constructed asymptotic SCB for
the mean function of
the functional regression model by using piecewise constant
spline estimation and a strong approximation result.
For functional data, Degras~\cite{Degras2011} constructed asymptotic
SCB for the mean function of
the functional linear model without considering any covariate,
while
Zhang and Chen~\cite{Zhang2007} adopted the method of ``smoothing
first, then estimation''
and propose a global test statistic for testing $B_j$, but
their results cannot be used for constructing SCB for $B_j$.
Recently, Cardot et al.~\cite{Cardot2010}, Cardot and Josserand~\cite
{Cardot2011} built asymptotic SCB for Horvitz--Thompson estimators for
the mean function, but their models and estimation methods differ
significantly from ours.

In this paper, we propose an estimation procedure for the
multivariate varying coefficient model (\ref{Tract3eq1}) by using
local linear regression techniques, and derive a simultaneous
confidence band for the regression coefficient functions. We further
develop a test for linear hypotheses of coefficient functions. The
major aim of this paper is to investigate the theoretical properties
of the proposed estimation procedure and
test statistics. The theoretical development is challenging, but of
great interest for
carrying out statistical inferences on $B_j$.\vadjust{\goodbreak} The major
contributions of this paper are summarized as follows.
We first establish the weak convergence of the local linear estimator
of $B_j$, denoted by $\hat B_j$,
by using advanced empirical process methods \cite
{VaartWellner1996,Kosorok2008}. We further derive the bias and
asymptotic variance of
$\hat B_j$. These results provide insight into how the
direct estimation procedure for $B_j$ using observations from
all subjects outperforms the estimation procedure with the strategy
of ``smoothing first, then estimation.'' After calculating $\hat B_j$,
we reconstruct
all individual functions $\eta_{ij}$ and establish their uniform
convergence rates.
We derive uniform convergence rates of the proposed estimate for the covariance
matrix of $\eta_{ij}$ 
and its associated eigenvalue and eigenvector functions by using
related results in Li and Hsing~\cite{LiHsing2010}. Using the weak
convergence of the local linear estimator of~$B_j$, we further
establish the asymptotic distribution of a global test statistic for
linear hypotheses of the regression coefficient functions, and
construct an asymptotic SCB for each varying coefficient function.

The rest of this paper is organized as follows. In Section~\ref{sec2}, we
describe MVCM and its estimation procedure.
In Section~\ref{sec3}, we
propose a global test statistic for linear hypotheses of the
regression coefficient functions and construct an asymptotic
SCB for each coefficient function. In Section~\ref{sec4}, we discuss the
theoretical properties of estimation and inference procedures.
Two sets of simulation studies are presented in Section~\ref{sec5} with the
known ground truth to
examine the finite sample performance of the global test statistic
and SCB for each individual varying coefficient function. In Section~\ref{sec6},
we use MVCM to
investigate the development of white matter diffusivities along the
genu tract of the corpus callosum in a clinical study of neurodevelopment.

\section{Estimation procedure}\label{sec2}

Throughout this paper, we assume that
${\ve\varepsilon}_i(s)=(\varepsilon_{i 1}(s),\ldots, \varepsilon_{iJ}(s))^T$
and ${\ve\eta}_i(s)=(\eta_{i 1}(s),\ldots, \eta_{iJ}(s))^T$ are
mutually independent, and ${\ve\eta}_i(s)$ and ${\ve\varepsilon}_i(s)$
are independent and identical copies of SP$({\mathbf0}, \Sigma_{\eta})$
and SP$({\mathbf0}, \Sigma_\varepsilon)$, respectively, where SP$(\mu,
\Sigma)$ denotes a stochastic process vector with mean function
$\mu(t)$ and covariance function $\Sigma(s, t)$. Moreover,
$\ve\varepsilon_i(s)$ and $\ve\varepsilon_i(t)$ are assumed to be
independent for $s\not=t$, and $\Sigma_\varepsilon(s, t)$ takes the
form of $S_\varepsilon(t) {\mathbf1}(s=t)$, where $S_\varepsilon
(t)=(s_{\varepsilon, jj'}(t))$ is a $J\times J$ matrix of functions of
$t$ and ${\mathbf1}(\cdot)$
is an indicator function. Therefore, the covariance structure of
${\mathbf y}_i(s)$, denoted by $\Sigma_y(s, t)$, is given by
%
\begin{equation}
\label{Tract3eq2} \Sigma_y(s, t)=\operatorname{Cov}\bigl({\mathbf
y}_i(s), {\mathbf y}_i(t)\bigr)= \Sigma_{\eta}(s,
t)+S_\varepsilon(t){\mathbf1}(s=t).
\end{equation}

\subsection{Estimating varying coefficient functions}\label{sec2.1}

We employ local linear regression~\cite{Fan1996} to
estimate the coefficient functions $B_j$. Specifically, we apply
the Taylor expansion for
$B_j(s_m)$ at $s$ as follows:
%
\begin{equation}
\label{Tract3eq5} B_j(s_m) \approx B_j(s)+
\dot B_j(s) (s_m-s) =A_j(s){\mathbf
z}_{h_{1j}}(s_m-s),
\end{equation}
where ${\mathbf z}_{h}(s_m-s)=(1, (s_m-s)/h)^T$ and $A_j(s)=[B_j(s)
h_{1j} \dot B_j(s)]$ is a $p\times2$ matrix, in which $\dot
B_j(s)=(\dot b_{j1}(s),\ldots, \dot b_{jp}(s))^T$\vadjust{\goodbreak} is a $p\times1$
vector and $\dot b_{jl}(s)=d b_{jl}(s)/ds$ for $l=1,\ldots, p$. Let
$K(s)$\vspace*{1pt} be a kernel function and $K_h(s)=\break h^{-1} K(s/h)$ be the rescaled
kernel function with a bandwidth $h$.
We estimate $A_j(s)$ by minimizing the following weighted least
squares function:
%
\begin{equation}
\label{Tract3eq6} \sum_{i=1}^n\sum
_{m=1}^{M}\bigl[y_{ij}(s_m)-{
\mathbf x}_{i}^TA_j(s){\mathbf z}_{h_{1j}}(s_m-s)
\bigr]^2K_{h_{1j}}(s_m-s).
\end{equation}
%

Let us now introduce some matrix operators.
Let ${\mathbf a}^{\otimes2}={\mathbf a}{\mathbf a}^T$ for any vector ${\mathbf a}$
and $C\otimes D$ be the Kronecker product of two matrices $C$ and
$D$. For an $M_1\times M_2$ matrix $C=(c_{jl})$, denote
$\operatorname{vec}(C)=(c_{11},\ldots, c_{1 M_2},\ldots, c_{M_11},\ldots, c_{M_1 M_2})^T$. Let $\hat{A}_j(s)$ be the minimizer of
(\ref{Tract3eq6}). Then
%
\begin{equation}
\label{Tract3eq7} \operatorname{vec}\bigl(\hat{A}_j(s)\bigr)=
\Sigma(s,h_{1j})^{-1} \sum_{i=1}^n
\sum_{m=1}^{M} K_{h_{1j}}(s_m-s)
\bigl[{\mathbf x}_{i}\otimes{\mathbf z}_{h_{1j}}(s_m-s)
\bigr] y_{ij}(s_m),\hspace*{-35pt}
\end{equation}
where $\Sigma( s, h_{1j})=\sum_{i=1}^n\sum_{m=1}^{M}
K_{h_{1j}}(s_m-s)[{\mathbf x}_{i}^{\otimes2}\otimes{\mathbf
z}_{h_{1j}}(s_m-s)^{\otimes2}]$. Thus, we have
%
\begin{equation}
\label{Tract3eq8} \hat B_j(s)=\bigl(\hat b_{j1}(s),\ldots, \hat b_{jp}(s)\bigr)^T=\bigl[{ I}_p
\otimes(1, 0)\bigr] \operatorname{vec}\bigl(\hat{A}_j(s)\bigr),
\end{equation}
where ${I}_p$ is a $p\times p$ identity matrix.

In practice, we may select the bandwidth $h_{1j}$ by using
leave-one-curve-out cross-validation.
Specifically, for each $j$, we pool the data from
all $n$ subjects and select a bandwidth $h_{1j}$, denoted by $\hat
h_{1j}$, by minimizing the cross-validation score given by
%
\begin{equation}
\label{Tract3eq9} \operatorname{CV}(h_{1j})=(nM)^{-1}\sum
_{i=1}^n\sum_{m=1}^{M}
\bigl[y_{ij}(s_m)- {\mathbf x}_i^T
\hat B_j(s_m, h_{1j})^{(-i)}
\bigr]^2,
\end{equation}
where $\hat B_j(s, h_{1j})^{(-i)}$ is the local linear estimator of
$B_j(s)$ with the bandwidth $h_{1j}$ based on data excluding
all the observations from the $i$th subject.

\subsection{Smoothing individual functions}\label{sec2.2}

By assuming certain smoothness conditions on $\eta_{ij}(s)$, we also
employ the local linear regression technique to estimate
all individual functions $\eta_{ij}(s)$
\cite{Fan1996,Wand1995,Wu2006,Ramsay2005,Welsh2006,Zhang2007}.
Specifically, we have the Taylor expansion for ${\eta}_{ij}(s_m)$
at $s$,
%
\begin{equation}
\label{Tract3eq10} \eta_{ij}(s_m)\approx{\mathbf
d}_{ij}(s)^T{\mathbf z}_{h_{2j}}(s_m-s),
\end{equation}
where ${\mathbf d}_{ij}(s)=({\eta}_{ij}(s),
h_{2j}\dot{\eta}_{ij}(s))^T$ is a $ 2\times1$ vector. We
develop an algorithm to estimate ${\mathbf d}_{ij}(s)$ as follows.
For each $i$ and $j$, we estimate ${\mathbf d}_{ij}(s)$ by minimizing the
weighted least squares function:
%
\begin{equation}
\label{Tract3eq11} \sum_{m=1}^{M}
\bigl[{y}_{ij}(s_m)-{\mathbf x}_i^T
\hat B_{j}(s_m)-{\mathbf d}_{ij}(s)^T{
\mathbf z}_{h_{2j}}(s_m-s)\bigr]^2K_{h_{2j}}(s_m-s).
\end{equation}
Then, $\eta_{ij}(s)$ can be estimated by
%
\begin{eqnarray}
\label{Tract3eq12}
\hat\eta_{ij}(s)&=&(1,0)\hat{\mathbf d}_{ij}(s)\nonumber\\[-8pt]\\[-8pt]
&=&
\sum_{m=1}^{M} \tilde K_{h_{2j}}(s_{m}-s)
\bigl[{ y}_{ij}(s_m)-{\mathbf x}_i^T
\hat B_{j}(s_m)\bigr],\nonumber
\end{eqnarray}
where $\tilde K_{h_{2j}}(s)$ are the empirical
equivalent kernels and $\hat{\mathbf d}_{ij}(s)$ is given by
\begin{eqnarray*}
\hat{\mathbf d}_{ij}(s)&=&\Biggl[\sum_{m=1}^{M}K_{h_{2j}}(s_m-s){
\mathbf z}_{h_{2j}}(s_m-s)^{\otimes2}\Biggr]^{-1}
\\
&&{} \times\sum_{m=1}^{M}
K_{h_{2j}}(s_m-s){\mathbf z}_{h_{2j}}(s_m-s)
\bigl[{ y}_{ij}(s_m)-{\mathbf x}_i^T
\hat B_{j}(s_m)\bigr].
\end{eqnarray*}
Finally, let ${ S}_{ij}$ be the smoother matrix for the $j$th
measurement of the $i$th subject~\cite{Fan1996}, we can obtain
%
\begin{equation}
\label{Tract3eq13} \hat{\ve\eta}_{ij}=\bigl(\hat\eta_{i
j}(s_1),\ldots, \hat\eta_{ij}(s_{M})\bigr)^T={
S}_{ij} {R}_{i j},
\end{equation}
where $R_{i j}=(y_{ij}(s_1)-{\mathbf x}_i^T\hat B_{j}(s_1),\ldots,
y_{ij}(s_{M})-{\mathbf x}_i^T\hat B_{j}(s_{M}))^T$.

A simple and efficient way to obtain $h_{2j}$ is to use generalized
cross-validation method.
For each $j$, we pool the data from all $n$ subjects and select the
optimal bandwidth $h_{2j}$, denoted by $\hat h_{2j}$, by
minimizing the generalized cross-validation score given by
%
\begin{equation}
\label{Tract3eq14} \operatorname{GCV}(h_{2j})=\sum
_{i=1}^n \frac
{R_{ij}^T(I_{M}-S_{ij})^T(I_{M}-S_{ij})R_{ij}
}{[1-M^{-1}\operatorname{tr}(S_{ij})]^2}.
\end{equation}
Based on $\hat h_{2j}$, we can use (\ref{Tract3eq12}) to
estimate $\eta_{ij}(s)$ for all $i$ and $j$.

\subsection{Functional principal component analysis}\label{sec2.3}

We consider a spectral decomposition of $\Sigma_{\eta}(s, t)=(\Sigma_{\eta, jj'}(s, t))$ and its approximation.
According to Mercer's theorem~\cite{Mercer1909}, if $\Sigma_\eta(s,
t)$ is continuous on $[0, 1]\times[0, 1]$, then
$\Sigma_{\eta, jj}(s,
t)$ admits a spectral decomposition. Specifically,
we have
%
\begin{equation}
\label{Tract3eq3} \Sigma_{\eta, jj}(s, t)=\sum_{l=1}^\infty
\lambda_{j l} \psi_{j
l}(s)\psi_{j l}(t)
\end{equation}
for $j=1,\ldots, J$,
where $\lambda_{j 1}\geq\lambda_{j 2}\geq\cdots\geq0$ are
ordered values of the eigenvalues of a linear operator determined by
$\Sigma_{\eta, jj}$ with $\sum_{l=1}^\infty\lambda_{j l}<\infty$
and the $\psi_{jl}(t)$'s are the corresponding orthonormal
eigenfunctions (or principal components)
\cite{LiHsing2010,MR2212572,MR2278365}. 
The eigenfunctions form an orthonormal system on the space of
square-integrable functions on $[0, 1]$, and
$\eta_{ij}(t)$ admits the Karhunen--Loeve expansion as
$\eta_{ij}(t)=\sum_{l=1}^\infty\xi_{ijl}\psi_{j l}(t)$, where
$\xi_{ijl}=\int_0^{1} \eta_{ij}(s)\psi_{j l}(s)\,ds$ is referred to as
the $(jl)$th functional principal component scores of the $i$th
subject.
For each fixed $(i, j)$, the $\xi_{ijl}$'s are uncorrelated random
variables with $E(\xi_{ijl})=0$ and $E(\xi_{ijl}^2)=\lambda_{j l}$.
Furthermore, for $j\not=j'$, we have
\[
\Sigma_{\eta, jj'}(s,t)=\sum_{l=1}^\infty
\sum_{l'=1}^\infty E(\xi_{ij l}
\xi_{ij' l'}) \psi_{j l}(s)\psi_{j' l'}(t).
\]
%

After obtaining $\hat{\ve\eta}_{i}(s)=(\hat\eta_{i1}(s),\ldots,
\hat\eta_{iJ}(s))^T$, we estimate
$\Sigma_\eta(s, t)$ by using the empirical covariance of the
estimated $\hat{\ve\eta}_{i}(s)$ as follows:
\[
\hat\Sigma_\eta(s, t)=(n-p)^{-1}\sum
_{i=1}^n \hat{\ve\eta}_{i}(s) \hat{\ve
\eta}_{i}(t)^T.
\]
Following Rice and Silverman~\cite{MR1094283}, we can calculate the spectral
decomposition of $\hat\Sigma_{\eta, jj}(s, t)$ for each $j$ as
follows:
%
\begin{equation}
\label{Tract3eq15} \hat\Sigma_{\eta, jj}(s, t)=\sum
_l \hat\lambda_{j l} \hat\psi_{j
l}(s)
\hat\psi_{j l}(t),
\end{equation}
where $\hat\lambda_{j 1}\geq\hat\lambda_{j 2}\geq\cdots\geq0$
are estimated eigenvalues and the $\hat\psi_{j l}(t)$'s are the
corresponding estimated principal components. Furthermore, the $(j,
l)$th functional principal component scores can be computed using
$\hat\xi_{ijl}=\break \sum_{m=1}^{M}\hat\eta_{ij}(s_m)\*\hat\psi_{j
l}(s_m)(s_{m}-s_{m-1})$ for $i=1,\ldots, n$. We further show the
uniform convergence
rate of $\hat\Sigma_\eta(s, t)$ and its associated eigenvalues and
eigenfunctions. This result is useful for constructing the global
and local test statistics for testing the covariate effects.


\section{Inference procedure}\label{sec3}

In this section, we study global tests for linear hypotheses of
coefficient functions and SCB for each varying coefficient
function. They are essential for statistical inference on the
coefficient functions.\looseness=1

\subsection{Hypothesis test}\label{sec3.1}
Consider the linear hypotheses of ${\mathbf B}(s)$ as follows:
%
\begin{equation}
\label{Tract3eq18}
H_0\dvtx  {\mathbf C}\operatorname{vec}\bigl({\mathbf B}(s)\bigr)={\mathbf
b}_0(s)\qquad\mbox{for all } s \quad\mbox{vs.}\quad H_1\dvtx  {\mathbf C}
\operatorname{vec}\bigl({\mathbf B}(s)\bigr)\not={\mathbf b}_0(s),\hspace*{-35pt}
\end{equation}
where ${\mathbf B}(s)=[ B_{1}(s),\ldots, B_{J}(s)]$, ${\mathbf C}$ is a
$r\times Jp$ matrix with rank $r$ and ${\mathbf b}_0(s)$ is a given
$r\times
1$ vector of functions.
Define a global test statistic $S_n$ as
%
\begin{equation}
\label{Tract3eq21} S_n=\int_0^{1} {
\mathbf d}(s)^T \bigl[{\mathbf C}\bigl(\hat\Sigma_\eta(s, s)\otimes
\hat\Omega_X^{-1}\bigr){\mathbf C}^T
\bigr]^{-1} {\mathbf d}(s) \,ds,
\end{equation}
where $\hat\Omega_X=\sum_{i=1}^n {\mathbf x}_i^{\otimes2}$ and ${\mathbf
d}(s)= {\mathbf C}\operatorname{vec}(\hat{\mathbf B}(s)-\operatorname{bias}(\hat{\mathbf
B}(s)))-{\mathbf b}_0(s)$.



To calculate $S_n$, we need to estimate the bias of $\hat B_{j}(s)$
for all $j$. Based on (\ref{Tract3eq8}), we have
%
\begin{eqnarray}\label{Tract3eq20}
&&
\operatorname{bias}\bigl(\hat B_{j}(s)\bigr)\nonumber\\
&&\qquad= \bigl[{ I}_p
\otimes(1, 0)\bigr] \nonumber\\[-8pt]\\[-8pt]
&&\qquad\quad{}\times\operatorname{vec}\Biggl(\Sigma( s,h_{1j})^{-1}
\sum_{i=1}^n\sum
_{m=1}^{M} K_{h_{1j}}(s_m-s)\bigl[{\mathbf x}_{i}\otimes{\mathbf z}_{h_{1j}}(s_m-s)
\bigr]
\nonumber\\
&&\qquad\quad\hspace*{121pt}{} \times  {\mathbf x}_{i}(s_m)^T\bigl[B_j(s_m)-
\hat A_j(s){\mathbf z}_{h_{1j}}(s_m-s)\bigr]\Biggr).
\nonumber
\end{eqnarray}
By using Taylor's expansion, we have
\[
B_j(s_m)- \hat A_j(s){\mathbf
z}_{h_{1j}}(s_m-s)\approx2^{-1}\ddot
B_j(s) (s_m-s)^2+6^{-1} \dddot
B_j(s) (s_m-s)^3,
\]
where $\ddot B_j(s)=d^2 B_j(s)/ds^2$ and $\dddot B_j(s)=d^3
B_j(s)/ds^3$. Following the pre-asymptotic substitution method of
Fan and Gijbels~\cite{Fan1996}, we replace $B_j(s_m)- \hat A_j(s){\mathbf
z}_{h_{1j}}(s_m-s)$ by $ 2^{-1}\hhat{\ddot B}_j(s)(s_m-s)^2+6^{-1}
\hhat{\dddot B}_j(s) (s_m-s)^2$, in which $\hhat{\ddot
B}_j(s)$ and $\hhat{\dddot B}_j(s)$ are estimators obtained by
using local cubic fit with a pilot bandwidth selected in (\ref{Tract3eq9}).


It will be shown below that the asymptotic distribution of $S_n$ is
quite complicated,
and it is difficult to directly approximate the percentiles
of $S_n$ under the null hypothesis. Instead, we propose
using a wild bootstrap method to obtain critical values of $S_n$. The
wild bootstrap consists of the following three steps:
%
\begin{step}\label{step1}
Fit model (\ref{Tract3eq1}) under the null
hypothesis $H_0$, which yields $\hat B^*(s_m)$, $\hat
{\ve\eta}_{i,0}^*(s_m)$ and $\hat{\ve\varepsilon}_{i,0}^*(s_m)$ for
$i=1,\ldots, n$ and $m=1,\ldots, M$.
\end{step}
%
\begin{step}\label{step2}
Generate a random sample $\tau_i^{(g)}$ and $\tau_i(s_m)^{(g)}$ from a $N(0,
1)$ generator for $i=1,\ldots, n$ and $m=1,\ldots, M$ and then
construct
\[
\hat{\mathbf y}_i(s_m)^{(g)}=\hat
B^*(s)^T{\mathbf x}_i+\tau^{(g)}_i\hat
{\ve\eta}_{i,0}^*(s_m)+\tau_i(s_m)^{(g)}
\hat{\ve\varepsilon}_{i,0}^*(s_m).
\]
Then, based on $\hat{\mathbf y}_i(s_m)^{(g)}$, we recalculate $\hat{\mathbf
B}(s)^{(g)}$, $\operatorname{bias}(\hat{\mathbf
B}(s)^{(g)})$ and ${\mathbf d}(s)^{(g)}={\mathbf C}\operatorname{vec}(\hat{\mathbf
B}(s)^{(g)}-\operatorname{bias}(\hat{\mathbf B}(s)^{(g)}))-{\mathbf b}_0(s)$. We
also note that ${\mathbf C}\operatorname{vec}(\hat{\mathbf B}(s)^{(g)})\approx
{\mathbf b}_0$ and ${\mathbf C}\operatorname{vec}(\operatorname{bias}(\hat{\mathbf
B}(s)^{(g)}))\approx{\mathbf0}$. Thus, we can drop the term
$\operatorname{bias}(\hat{\mathbf B}(s)^{(g)})$ in ${\mathbf d}(s)^{(g)}$ for
computational efficiency. Subsequently, we compute
\[
S_n^{(g)}=n\int_0^{1} {
\mathbf d}(s)^{(g)T} \bigl[{\mathbf C}\bigl(\hat\Sigma_\eta (s, s)
\otimes\hat\Omega_X^{-1}\bigr){\mathbf C}^T
\bigr]^{-1} {\mathbf d}(s)^{(g)} \,ds.
\]
\end{step}
%
\begin{step}\label{step3}
Repeat Step~\ref{step2} $G$ times to obtain $\{ S_n^{(g)}\dvtx
g=1,\ldots, G\}$, and then calculate $ p= G^{-1}
\sum_{g=1}^G 1(S_n^{(g)}\geq S_n). $ If $p$ is smaller than
a pre-specified significance level $\alpha$, say 0.05, then one
rejects the null hypothesis $H_0$.
\end{step}

\subsection{Simultaneous confidence bands}\label{sec3.2}

Construction of SCB for coefficient functions is of great interest
in statistical inference for model (\ref{Tract3eq1}).
For a given confidence level $\alpha$, we construct SCB for each
$b_{jl}(s)$ as follows:
%
\begin{equation}
\label{Tract3eq22} P\bigl( \hat b^{L, \alpha}_{jl}(s)<b_{jl}(s)<
\hat b^{U,
\alpha}_{jl}(s) \mbox{ for all }s\in[0, 1]\bigr)=1-\alpha,
\end{equation}
where $\hat b^{L, \alpha}_{jl}(s)$ and $\hat b^{U, \alpha}_{jl}(s)$
are the lower and upper limits of SCB.
Specifically, it will be shown below that
a $1-\alpha$ simultaneous confidence band for $b_{jl}(s)$ is given as follows:
%
\begin{equation}
\label{Tract3eq25}\quad \biggl( \hat b_{jl}(s)-\operatorname{bias}\bigl(\hat
b_{jl}(s)\bigr)-\frac
{C_{jl}(\alpha)}{\sqrt{n}}, \hat b_{jl}(s)-\operatorname{bias}
\bigl(\hat b_{jl}(s)\bigr)+\frac{C_{jl}(\alpha
)}{\sqrt{n}} \biggr),
\end{equation}
where $C_{jl}(\alpha)$ is a scalar.
Since the calculation of $\hat b_{jl}(s)$ and $\operatorname{bias}(\hat
b_{jl}(s))$ has been discussed in
(\ref{Tract3eq8}) and (\ref{Tract3eq20}),
the next issue is to determine $C_{jl}(\alpha)$.

Although there are several methods of determining $C_{jl}(\alpha)$
including random field theory~\cite{Worsley2004,SunLoader1994}, we
develop an efficient resampling method to approximate
$C_{jl}(\alpha)$ as follows~\cite{Zhu2007a,Kosorok2003}:
\begin{itemize}
\item We calculate $\hat{r}_{ij}(s_m)=y_{ij}(s_m)-{\mathbf
x}_i^T\hat B_{j}(s_m)$ for all $i, j$, and $m$.
\item For $g=1,\ldots, G$, we independently
simulate $\{\tau_i^{(g)}\dvtx  i=1,\ldots, n\}$ from $N(0, 1)$ and
calculate a stochastic process $G_{j}(s)^{(g)}$ given by
\begin{eqnarray*}
\hspace*{-5pt}
&&
\sqrt{n} \bigl[{ I}_p\otimes(1, 0)\bigr]\\
\hspace*{-5pt}
&&\quad{}\times \operatorname{vec}
\Biggl(\Sigma(s, h_{1j})^{-1} \sum
_{i=1}^n\tau_i^{(g)}\sum
_{m=1}^{M} K_{h_{1j}}(s_m-s)
\bigl[{\mathbf x}_{i}\otimes{\mathbf z}_{h_{1j}}(s_m-s)
\bigr] \hat{r}_{ij}(s_m)\Biggr).
\end{eqnarray*}
%
%
\item We calculate $\sup_{s\in[0, 1]}|{\mathbf
e}_lG_{j}(s)^{(g)}|$ for all $g$, where ${\mathbf e}_l$ be a $p\times1$
vector with the $l$th element 1 and 0
otherwise, and use their $1-\alpha$ empirical
percentile to estimate $C_{jl}(\alpha)$.
\end{itemize}

\section{Asymptotic properties}\label{sec4}

In this section, we systematically examine the asymptotic properties of
$\hat{\mathbf B}(s)$, $\hat\eta_{ij}(s)$,
$\hat\Sigma_\eta(s, t)$ and $S_n$ developed in Sections~\ref{sec2} and~\ref{sec3}.
Let us first define some notation. Let $u_r(K)=\int t^r K(t)\,dt$ and
$v_r(K)=\int t^r K^2(t) \,dt$, where $r$ is any integer. For
any\vspace*{1pt} smooth functions $f(s)$ and $g(s, t)$, define $\dot
f(s)=df(s)/ds$, $\ddot f(s)=d^2f(s)/ds^2$, $\dddot f(s)=d^3f(s)/ds^3$
and $g^{(a,b)}(s, t)=\partial^{a+b} g(s, t)/\partial^{a}
s\,\partial^{b} t$, where $a$ and $b$ are any nonnegative integers. Let
${\mathbf H}=\operatorname{diag}( h_{11},\ldots, h_{1J})$, ${\mathbf
B}(s)=[ B_{1}(s),\ldots, B_{J}(s)]$, $\hat{\mathbf B}(s)=[\hat
B_{1}(s),\ldots, \break\hat B_{J}(s)]$ and $\ddot{\mathbf B}(s)=[\ddot
B_{1}(s),\ldots, \ddot B_{J}(s)]$, where $\ddot B_j(s)=(\ddot
b_{j1}(s),\ldots, \ddot b_{jp}(s))^T$ is a $p\times1$ vector. Let
${\mathcal S}=\{s_1,\ldots, s_M\}$.

\subsection{Assumptions}\label{sec4.1}

Throughout the paper, the following assumptions are needed to
facilitate the technical details, although they may not be the
weakest conditions.
We need to introduce some notation. Let $N(\mu, \Sigma)$ be a normal
random vector with mean $\mu$ and covariance
$\Sigma$. Let $
\Omega_1(h, s)=\int (1, h^{-1}(u-s))^{\otimes2}K_h(u-s)\pi(u)\,du.
$
Moreover, 
we do not distinguish the
differentiation and continuation at the boundary points from those
in the interior of $[0, 1]$. For instance, a continuous function
at the boundary of $[0, 1]$ means that this function is left
continuous at $0$ and right continuous at $1$.


\renewcommand{\theassumption}{(C\arabic{assumption})}
\begin{assumption}\label{assumC1}
For all $j=1,\ldots, J$,
$\sup_{s_m} E[|\varepsilon_{ij}(s_m)|^{q}]<\infty$ for some $q>4$ and
all grid points $s_m$.
\end{assumption}

\begin{assumption}\label{assumC2}
Each component of $\{ \ve\eta(s)\dvtx  s\in[0, 1]\}$,
$\{ \ve\eta(s)\ve\eta(t)^T\dvtx \break (s, t)\in[0, 1]^2\}$
and $\{ {\mathbf x}\ve\eta^T(s)\dvtx  s\in[0, 1]\}$ are Donsker classes.
\end{assumption}

\begin{assumption}\label{assumC3}
The covariate\vspace*{1pt} vectors ${\mathbf x}_i$'s are
independently and identically
distributed with $E{\mathbf x}_i=\mu_x$ and $\|{\mathbf x}_{i}\|_\infty
<\infty$. Assume that $E[{\mathbf x}_{i}^{\otimes2}]
=\Omega_X$ is positive definite.
\end{assumption}

\begin{assumption}\label{assumC4}
The grid points ${\mathcal S}=\{s_m, m=1,\ldots,
M\}$ are randomly generated from a density function $\pi(s)$. Moreover,
$\pi(s)>0$ for all $s\in[0, 1]$ and $\pi(s)$ has continuous
second-order derivative with the bounded support $[0, 1]$.
\end{assumption}

{\renewcommand{\theassumption}{(C4\textup{b})}
\begin{assumption}\label{assumC4b}
The grid\vspace*{1pt} points ${\mathcal S}=\{s_m, m=1,\ldots, M\}$ are prefixed according to $\pi(s)$ such that $\int_0^{s_m}\pi(s)\,ds=m/M$ for $M\geq m\geq1$. Moreover, $\pi(s)>0$ for
all $s\in[0, 1]$ and $\pi(s)$ has continuous second-order derivative
with the bounded support $[0, 1]$.
\end{assumption}}

\setcounter{assumption}{4}
\begin{assumption}\label{assumC5}
The kernel function $K(t)$ is a symmetric
density function with a compact
support $[-1, 1]$, and is Lipschitz continuous.
Moreover,
$0<\inf_{h\in(0, h_0], s\in[0, 1]} \operatorname{det}(\Omega_1(h, s))$,
where $h_0>0$ is a small scalar
and
$\operatorname{det} (\Omega_1(h, s))$ denotes the determinant of $\Omega_1(h, s)$.
\end{assumption}

\begin{assumption}\label{assumC6}
All components of ${\mathbf B}(s)$ have continuous
second derivatives on $[0, 1]$.
\end{assumption}

\begin{assumption}\label{assumC7}
Both $n$ and $M$ converge to $\infty$, $\max_{j}h_{1j}=o(1)$,
$Mh_{1j}\rightarrow\infty$ and
$\max_j h_{1j}^{-1}|{\log h_{1j}}|^{1-2/q_1}\leq M^{1-2/q_1}$ for $j=1,\ldots, J$, where $q_1\in(2, 4)$.
\end{assumption}

{\renewcommand{\theassumption}{(C7\textup{b})}
\begin{assumption}\label{assumC7b}
Both $n$ and $M$ converge to $\infty$, $\max_{j}h_{1j}=o(1)$,
$Mh_{1j}\rightarrow\infty$ and $\log(M)=o(Mh_{1j})$. There exists a
sequence of $\gamma_n>0$ such that $\gamma_n\rightarrow\infty$,
$\max_j n^{1/2}\gamma_n^{1-q}h_{1j}^{-1}=o(1)$ and $n^{-1/2}\gamma_n\log(M)=o(1)$.
\end{assumption}}

\setcounter{assumption}{7}
\begin{assumption}\label{assumC8}
For all $j$, $\max_j (h_{2j})^{-4}(\log n/n)^{1-2/q_2}=o(1)$ for
$q_2\in(2, \infty)$,
$\max_{j}h_{2j}=o(1)$, and $Mh_{2j}\rightarrow\infty$
for $j=1,\ldots, J$.
\end{assumption}



\begin{assumption}\label{assumC9a}
The sample path
of $\eta_{ij}(s)$ has continuous second-or\-der derivative on $[0, 1]$
and
$E[\sup_{s\in[0, 1]}\|{\ve\eta}(s)\|_2^{r_1}]\!<\!\infty$ and $E\{\sup_{s\in[0, 1]}[\|\dot{\ve\eta}(s)\|_2\!+\!\|\ddot{\ve\eta
}(s)\|_2]^{r_2}\}<\infty$ for some $r_1, r_2\in(2, \infty)$, where
$\|\cdot\|_2$ is the Euclidean norm.
\end{assumption}

{\renewcommand{\theassumption}{(C9\textup{b})}
\begin{assumption}\label{assumC9b}
$E[\sup_{s\in[0, 1]}\|{\ve\eta}(s)\|_2^{r_1}]<\infty$ for some
$r_1\in(2, \infty)$ and\vspace*{1pt}
all components of $\Sigma_\eta(s, t)$ have continuous second-order
partial derivatives with respect to $(s, t)\in[0, 1]^2$ and
$\inf_{s\in[0, 1]}\Sigma_\eta(s, s)>0$.
\end{assumption}}




\setcounter{assumption}{9}
\begin{assumption}\label{assumC10}
There is a positive fixed integer $E_j<\infty$
such that $\lambda_{j, 1}>\cdots>\lambda_{j, E_j}>\lambda_{j,
E_j+1}\geq\cdots\geq0$ for $j=1, \ldots, J$.
\end{assumption}



\begin{Remark*}
Assumption~\ref{assumC1} requires the uniform bound on the
high-order moment of $\varepsilon_{ij}(s_m)$ for all grid points $s_m$.
Assumption~\ref{assumC2} avoids smoothness conditions on the sample path $\ve
\eta(s)$, which are commonly assumed
in the literature
\cite{Degras2011,Zhang2007,MR2278365}.
Assumption~\ref{assumC3} is a relatively weak condition on the covariate vector,
and the boundedness of $\|{\mathbf x}_i\|_2$ is not essential.
Assumption~\ref{assumC4} is a weak condition on the random grid points.
In many neuroimaging applications, $M$ is often much larger than $n$ and
for such large $M$, a regular grid of voxels
is fairly well approximated by voxels generated by a uniform
distribution in a compact subset of Euclidean space. For notational simplicity,
we only state the theoretical results for the random grid points
throughout the paper.
Assumption~\ref{assumC4b} is a weak condition on the fixed grid points.
We will prove several key results for the fixed grid point case in
Lemma 8 of the supplemental article~\cite{ZLK2012}.
The bounded support restriction on $K(\cdot)$
in Assumption~\ref{assumC5} is not essential and can be removed if we put
a restriction on the tail of $K(\cdot)$.
Assumption~\ref{assumC6} is the standard smoothness condition on ${\mathbf B}(s)$ in
the literature
\mbox{\cite
{MR1742497,MR1769751,MR1959093,MR2425354,MR2504204,MR1666699,Ramsay2005,MR2087972,MR1888349,Zhang2007,MR2523900}}.
Assumptions~\ref{assumC7} and~\ref{assumC8} on
bandwidths are similar to the conditions used in \cite
{LiHsing2010,EinmahlMason2000}.
Assumption~\ref{assumC7b} is a weak condition on $n$, $M$, $h_{1j}$ and $\gamma_n$ for the fixed grid point case. For instance,
if we set $\gamma_n=n^{1/2}\log(M)^{-1-c_0}$ for a positive scalar
$c_0>0$, then we have
$n^{1/2}\gamma_n^{1-q}h_{1j}^{-1}= n^{1-q/2}\log
(M)^{(1+c_0)(q-1)}h_{1j}^{-1}=o(1)$ and $n^{-1/2}\gamma_n\log(M)=
\log(M)^{-c_0}=o(1)$.
As shown in Theorem~\ref{theo1} below, if $h_{1j}
=O((nM)^{-1/5})$ and $\gamma_n=n^{1/2}\log(M)^{-1-c_0}$,
$n^{1/2}\gamma_n^{1-q}h_{1j}^{-1}$ reduces to $n^{6/5-q/2}\log
(M)^{(1+c_0)(q-1)} M^{1/5}$.
For relatively large $q$ in Assumption~\ref{assumC1}, $n^{6/5-q/2}\log
(M)^{(1+c_0)(q-1)} M^{1/5}$ can converge to zero.
Assumptions~\ref{assumC9a} and~\ref{assumC3} are sufficient conditions of Assumption~\ref{assumC2}.
Assumption~\ref{assumC9b} on the sample path is the same as Condition C6 used in
\cite{LiHsing2010}.
Particularly, if we use the method for estimating $\Sigma_\eta(s,
s')$ considered
in Li and Hsing~\cite{LiHsing2010}, then the differentiability of $\ve
\eta(s)$ in Assumption~\ref{assumC9a} can be
dropped.
Assumption~\ref{assumC10} on
simple multiplicity of the first
$E_j$ eigenvalues
is only needed to investigate the asymptotic properties of eigenfunctions.
\end{Remark*}

\subsection{\texorpdfstring{Asymptotic properties of $\hat{\mathbf B}(s)$}
{Asymptotic properties of B(s)}}\label{sec4.2}

The following theorem establishes the weak
convergence of $\{\hat{B}(s), s\in[0, 1]\}$, which is essential
for constructing global test statistics and SCB for ${\mathbf
B}(s)$.
%
\begin{theorem}\label{theo1}
Suppose that Assumptions~\ref{assumC1}--\ref{assumC7}
hold. The following results hold:

\begin{longlist}
\item
$\sqrt{n}\{\operatorname{vec}(\hat{\mathbf B}(s)-{\mathbf B}(s) -0.5
\ddot{\mathbf B}(s){\mathbf U}_2(K; s, {\mathbf H}) {\mathbf H}^{2}[1+o_p(1)] )\dvtx  s\in
[0, 1]\}$ converges weakly to a centered Gaussian process
$G(\cdot)$ with covariance matrix $ \Sigma_\eta(s, s')\otimes
\Omega_X^{-1}$, where $\Omega_X=E[{\mathbf x}^{\otimes2}]$ and ${\mathbf
U}_2(K; s, {\mathbf H}) $ is a $J\times J$ diagonal matrix, whose
diagonal elements will be defined in Lemma~\ref{lemma5} in the \hyperref[app]{Appendix}.

\item The asymptotic bias and conditional variance 
of $\hat{B}_j(s)$ given ${\mathcal S}$ for
$s\in(0, 1)$ are given by $
0.5 h_{1j}^{2}u_2(K) \ddot B_j(s)[1+o_p(1)]
$
and
$ n^{-1}
\Sigma_{\eta, jj}(s, s) \Omega_X^{-1}[1+o_p(1)]$,
respectively.
\end{longlist}
\end{theorem}

\begin{Remarks*} (1) The major
challenge in proving Theorem~\ref{theo1}(i) is dealing with within-subject
dependence. This is because the dependence between
${\ve\eta}(s)$ and ${\ve\eta}(s')$ in the newly proposed
multivariate varying coefficient model does not converge to zero due
to the within-curve dependence. It is worth noting that for any
given $s$, the corresponding asymptotic normality of $\hat{\mathbf
B}(s)$ may be established by using related techniques in
Zhang and Chen~\cite{Zhang2007}.
%
%
However, the marginal asymptotic normality does not imply the weak
convergence of $\hat{\mathbf B}(s)$ as a stochastic process in $[0,
1]$, since we need to verify the asymptotic continuity of $\{\hat
{\mathbf B}(s)\dvtx  s\in[0, 1]\}$ to establish its weak convergence.
In addition, Zhang and Chen~\cite{Zhang2007}
considered ``smoothing first, then estimation,'' which requires a
stringent assumption such that $n = O(M^{4/5})$. Readers are
referred to Condition A.4 and Theorem 4 in Zhang and Chen~\cite
{Zhang2007} for
more details.
In contrast, directly estimating ${\mathbf B}(s)$ using local kernel
smoothing avoids such stringent assumption on the numbers of grid
points and subjects.


(2) Theorem~\ref{theo1}(ii) only provides us the asymptotic bias and conditional
variance of $\hat{B}_j(s)$ given ${\mathcal S}$ for the interior
points of $(0, 1)$. The asymptotic bias and conditional variance
at the boundary points $0$ and $1$ are given in Lemma~\ref{lemma5}.
The asymptotic bias of $\hat B_j(s)$ is of the order $h_{1j}^2$, as
the\vspace*{-1pt} one in nonparametric regression setting. Moreover, the
asymptotic conditional variance of $\hat B_j(s)$ has a complicated
form due to the within-curve dependence. The leading term in the
asymptotic conditional variance is of order $n^{-1}$, which
is slower than the standard nonparametric rate $(nM h_{1j})^{-1}$
with the assumption $h_{1j}\rightarrow0$ and $M h_{1j}\rightarrow
\infty$.

(3) Choosing an optimal bandwidth $h_{1j}$ is not a trivial task for
model (\ref{Tract3eq1}). Generally, any bandwidth $h_{1j}$ satisfying
the assumptions $h_{1j}\rightarrow0$ and $M h_{1j}\rightarrow\infty$
can ensure the weak convergence of $\{ \hat{\mathbf B}(s)\dvtx  s\in[0,
1]\}$. Based on the asymptotic bias and conditional variance of
$\hat{\mathbf B}(s)$, we can calculate an optimal bandwidth for estimating
${\mathbf B}(s)$, $h_{1j}
=O_p((nM)^{-1/5})$. In this case, $n^{-1}h_{1j}^2 $ and $(nM)^{-1}h_{1j}$
reduce to $O_p(n^{-7/5}M^{-2/5})$ and $(nM)^{-6/5}$, respectively,
and their contributions depend on the relative size of $n$ over
$M$.\vspace*{-2pt}
\end{Remarks*}

\subsection{\texorpdfstring{Asymptotic properties of $\hat\eta_{ij}(s)$}
{Asymptotic properties of eta ij(s)}}\label{sec4.3}

We next study the asymptotic bias and covariance of
$\hat\eta_{ij}(s)$ as follows. We distinguish between two cases. The first
one is conditioning on the design points in ${\mathcal S}$, ${\mathbf
X}$, and ${\ve\eta}$. The other is conditioning on the design
points in ${\mathcal S}$ and ${\mathbf X}$. We define $ K^{*}((s-t)/h)=
\int K(u)K(u+(s-t)/h)\,du. $\vspace*{-2pt}
%
\begin{theorem}\label{theo2}
Under Assumptions~\ref{assumC1} and~\ref{assumC3}--\ref{assumC8},
the following results hold for all $s\in(0,
L)$:
\begin{longlist}[(a)]
\item[(a)] Conditioning on $({\mathcal S}, {\mathbf X}, {\ve\eta})$, we have
\begin{eqnarray*}
&&\operatorname{Bias}\bigl[\hat\eta_{ij}(s)|{\mathcal S}, {\ve\eta}, {\mathbf
x}_i\bigr]
\\[-1pt]
&&\qquad=0.5 u_2(K)\bigl[\ddot\eta_{ij}(s)h_{2j}^{2}+{
\mathbf x}_i^T\ddot B_j(s_m)
h_{1j}^{2}\bigr] \bigl[1+o_p(1)
\bigr]+O_p\bigl(n^{-1/2}\bigr),
\\[-1pt]
&&\operatorname{Cov}\bigl[\hat\eta_{ij}(s), \hat\eta_{ij}(t)|{
\mathcal S}, {\ve \eta}, {\mathbf x}_i\bigr]
\\[-1pt]
&&\qquad=K^*\bigl((s-t)/h_{2j}\bigr)\pi(t)^{-1}(Mh_{2j})^{-1}O_p(1)-
{\mathbf x}_i^T\Omega_X^{-1}{\mathbf
x}_i(nM h_{1j} )^{-1}O_p(1).
\end{eqnarray*}

\item[(b)] The asymptotic bias and covariance of $\hat\eta_{ij}(s)$
conditioning on ${\mathcal S}$ and ${\mathbf X}$ are given by
\begin{eqnarray*}
&&\operatorname{Bias}\bigl[\hat\eta_{ij}(s)|{\mathcal S}, {\mathbf X}\bigr]=0.5
u_2(K) {\mathbf x}_i^T\ddot B_j(s_m)
h_{1j}^{2}\bigl[1+o_p(1)\bigr],\\[-25pt]
\end{eqnarray*}
\begin{eqnarray*}
&&\operatorname{Cov}\bigl(\hat\eta_{ij}(s)-\eta_{ij}(s), \hat
\eta_{ij}(t)-\eta_{ij}(t)|{\mathcal S}, {\mathbf X}\bigr)
\\[-1pt]
&&\qquad=\bigl[1+o_p(1)\bigr] \bigl[0.25 u_2(K)^2h_{2j}^{4}
\Sigma_{\eta, jj}^{(2, 2)}(s, t)\\[-1pt]
&&\qquad\quad\hspace*{53.5pt}{}
+K^*\bigl((s-t)/h_{2j}\bigr)
\pi(t)^{-1}(Mh_{2j})^{-1}O_p(1)\\[-1pt]
&&\hspace*{128.5pt}\qquad\quad{}+
n^{-1}{\mathbf x}_i^T\Omega_X^{-1}{
\mathbf x}_i\Sigma_{\eta, jj}(s, t)\bigr].
\end{eqnarray*}

\item[(c)] The mean integrated squared error (MISE) of all $\hat\eta_{ij}(s)$ is given by
%
\begin{eqnarray}\label{theorem2Eq1}
&& n^{-1}\sum_{i=1}^n \int
_0^{1}E\bigl\{\bigl[\hat\eta_{ij}(s)-
\eta_{ij}(s)\bigr]^2|{\mathcal S}\bigr\} \pi(s)\,ds
\nonumber\\[-1pt]
&&\qquad= \bigl[1+o_p(1)\bigr]\nonumber\\[-9pt]\\[-9pt]
&&\qquad\quad{}\times\biggl\{O\bigl((Mh_{2j})^{-1}
\bigr) + n^{-1}\int_0^{1}
\Sigma_{\eta, jj}(s, s)\pi(s) \,ds
\nonumber\\[-1pt]
&&\hspace*{17.2pt}\qquad\quad{} + 0.25 u_2^2(K)\int_0^{1}
\bigl[ \ddot B_j(s)^T\Omega_X\ddot
B_j(s) h_{1j}^{4} +\Sigma_{\eta, jj}^{(2, 2)}(s,
s)h_{2j}^{4}\bigr]\pi(s) \,ds \biggr\}.
\nonumber\vadjust{\goodbreak}
\end{eqnarray}

\item[(d)] The optimal bandwidth for minimizing MISE (\ref{theorem2Eq1}) is given by
%
\begin{equation}
\label{theorem2Eq2}
\hat h_{2j} =O\bigl(M^{-1/5}\bigr).
\end{equation}

\item[(e)] The first order LPK reconstructions $\hat\eta_{ij}(s)$
using $\hat h_{2j} $ in (\ref{theorem2Eq2}) satisfy
%
\begin{equation}
\label{theorem2Eq3} \sup_{s\in[0, 1]}\bigl|\hat \eta_{ij}(s)-
\eta_{ij}(s)\bigr|=O_p\bigl({\bigl|\log(M)\bigr|}^{1/2}M^{-2/5}+
h_{1j}^2+n^{-1/2}\bigr)
\end{equation}
for $i=1, \ldots, n$.
\end{longlist}
\end{theorem}

\begin{Remark*}
Theorem~\ref{theo2} characterizes the statistical properties of
smoothing individual curves ${\eta}_{ij}(s)$ after first estimating
$ B_j(s)$. Conditioning on individual curves $\eta_{i j}(s)$,
Theorem~\ref{theo2}(a) shows that $\operatorname{Bias}[\hat\eta_{ij}(s)|{\mathcal S},
{\mathbf X},
{\ve\eta}]$ is associated with $0.5 u_2(K) {\mathbf x}_i^T\ddot
B_j(s_m) h_{1j}^{2}$, which is the bias term of $\hat{B}_j(s)$
introduced in the estimation step, and $0.5
u_2(K)\ddot\eta_{ij}(s)h_{2j}^{2}$ is introduced in the smoothing
individual functions step. Without conditioning on $\eta_{i j}(s)$,
Theorem~\ref{theo2}(b) shows that the bias of $\hat\eta_{i j}(s)$ is mainly
controlled by the bias in the estimation step.
The MISE of $\hat\eta_{ij}(s)$ in Theorem~\ref{theo2}(c) is the sum of
$O_p(n^{-1}+ h_{1j}^{4})$ introduced by the estimation of $B_j(s)$ and
$O_p ((Mh_{2j})^{-1}+h_{2j}^{4})$ introduced by the
reconstruction of $\eta_{ij}(s)$.
The optimal bandwidth for minimizing the MISE of $\hat\eta_{i j}(s)$
is a standard bandwidth for LPK.
If we use the optimal bandwidth in Theorem~\ref{theo2}(d), then the MISE of
$\hat\eta_{i j}(s)$ can achieve the order of $n^{-1} + h_{1j}^{4}
+ M^{-4/5}$.
\end{Remark*}

\subsection{\texorpdfstring{Asymptotic properties of $\hat\Sigma_\eta(s,t)$}
{Asymptotic properties of Sigma eta(s,t)}}\label{sec4.4} In
this section, we study the asymptotic properties of $\hat\Sigma_\eta
(s,t)$ and its spectrum decomposition.
%
\begin{theorem}\label{theo3}
\textup{(i)} Under Assumptions~\ref{assumC1} and
\ref{assumC3}--\ref{assumC9a}, it follows that
\[
\sup_{(s, t)\in[0, 1]^2}\bigl|\hat\Sigma_\eta(s,t)-\Sigma_\eta(s,
t)\bigr|= O_p\bigl((Mh_{2j})^{-1}+
h_{1j}^2+h_{2j}^{2}+(\log
n/n)^{1/2}\bigr).\hspace*{-20pt}
\]

\begin{enumerate}[(ii)]
\item[(ii)]
Under Assumptions~\ref{assumC1} and~\ref{assumC3}--\ref{assumC10}, if the optimal
bandwidths $h_{mj} $ for $m=1, 2$ are used to reconstruct
$\hat B_{j}(s)$ and $\hat\eta_{ij}(s)$ for all $j$,
then for $l=1,\ldots, E_j$, we have the following results:
\begin{enumerate}[(a)]
\item[(a)]
$\int_0^{1}[\hat\psi_{jl}(s)-\psi_{jl}(s)]^2\,ds=O_p((Mh_{2j})^{-1}+ h_{1j}^2+h_{2j}^{2}+(\log n/n)^{1/2})$;
\item[(b)]
$ |\hat\lambda_{jl}-\lambda_{jl}|= O_p((Mh_{2j})^{-1}+ h_{1j}^2+h_{2j}^{2}+(\log n/n)^{1/2})$.
\end{enumerate}
\end{enumerate}
\end{theorem}

\begin{Remark*}
Theorem~\ref{theo3} characterizes the uniform weak convergence rates of
$\hat\Sigma_\eta(s, t)$, $\hat\psi_{j l}$ and $\hat\lambda_{j l}$ for
all $j$. It can be regarded as an extension of Theorems \mbox{3.3--3.6} in Li
and Hsing~\cite{LiHsing2010}, which established the uniform strong
convergence rates of these estimates with the sole presence of
intercept and $J=1$ in model~(\ref{Tract3eq1}). Another difference is
that Li and Hsing~\cite{LiHsing2010} employed all cross products
$y_{ij}y_{ik}$ for $j\not=k$ and then used the local polynomial kernel
to estimate $\Sigma_\eta(s, t)$. As discussed in Li and Hsing \cite
{LiHsing2010}, their approach can relax the assumption on the
differentiability of the individual curves. In contrast, following
Hall, M{\"u}ller and Wang~\cite{MR2278365} and Zhang and Chen
\cite{Zhang2007}, we directly fit a smooth curve to $\eta_{ij}(s)$ for
each $i$ and estimate $\Sigma_\eta(s, t)$ by the sample covariance
functions. Our approach is computationally simple and can ensure that
all $\hat\Sigma_{\eta, jj}(s, t)$ are positive semi-definite, whereas
the approach in Li and Hsing~\cite{LiHsing2010} cannot. This is
extremely important for high-dimensional neuroimaging data, which
usually contains a large number of locations (called voxels) on a
two-dimensional (2D) surface or in a 3D volume. For instance, the
number of $M$ can number in the tens of thousands to millions, and thus
it can be numerically infeasible to directly operate on
$\hat\Sigma_\eta(s, s')$.

We use $\tilde\Sigma_\eta(s, s')$ to denote the local linear
estimator of $\Sigma_\eta(s, s')$ proposed in
Li and Hsing~\cite{LiHsing2010}. Following the arguments in Li and
Hsing~\cite{LiHsing2010}, we can easily obtain the following result.
\end{Remark*}
%
\begin{corollary}\label{corollary1}
Under Assumptions~\ref{assumC1}--\ref{assumC8} and
\ref{assumC9b}, it follows that
\[
\sup_{(s, t)\in[0, 1]^2}\bigl|\tilde\Sigma_\eta(s,t)-\Sigma_\eta(s,
t)\bigr|= O_p\bigl( h_{1j}^2+h_{2j}^{2}+(
\log n/n)^{1/2}\bigr).
\]
\end{corollary}




\subsection{Asymptotic properties of the inference procedures}\label{sec4.5}

In this section, we discuss the asymptotic properties of the global
statistic $S_n$ and the critical values of SCB.
Theorem~\ref{theo1}
allows us to construct SCB for coefficient functions $b_{jl}(s)$.
It follows from Theorem~\ref{theo1} that
%
\begin{equation}
\label{Tract3eq23} \sqrt{n}\bigl[\hat b_{jl}(s)-b_{jl}(s)-
\operatorname{Bias}\bigl(\hat b_{jl}(s)\bigr)\bigr] \Rightarrow
G_{jl}(s),
\end{equation}
where $\Rightarrow$ denotes weak convergence of a sequence
of stochastic processes, and $G_{jl}(s)$ is a centered Gaussian process
indexed by
$s\in[0, 1]$. Therefore, let $X_{\mathbf C}(s)$ be a centered Gaussian
process, and we have
%
\begin{eqnarray}
\label{Tract3eq24} \bigl[{\mathbf C}\bigl(\hat\Sigma_\eta(s, s)\otimes
\hat\Omega_X^{-1}\bigr){\mathbf C}^T
\bigr]^{-1/2} {\mathbf d}(s)&\Rightarrow& X_{\mathbf C}(s),
\nonumber\\[-8pt]\\[-8pt]
\sup_{s\in[0, 1]}\bigl|\sqrt{n}\bigl[\hat b_{jl}(s)-b_{jl}(s)-
\operatorname{Bias}\bigl(\hat b_{jl}(s)\bigr)\bigr] \bigr|&\Rightarrow&
\sup_{s\in[0, 1]}\bigl|G_{jl}(s)\bigr|.\nonumber
\end{eqnarray}
We define $C_{jl}(\alpha)$ such that $P( \sup_{s\in[0,
1]}|G_{jl}(s)|\leq C_{jl}(\alpha))=1-\alpha$. Thus,
the confidence band given in (\ref{Tract3eq25}) is a $1-\alpha$
simultaneous confidence band for $b_{jl}(s)$.
%
\begin{theorem}\label{theo4}
If Assumptions~\ref{assumC1}--\ref{assumC9a} are true,
then we have
%
\begin{equation}
\label{Tract3NEW1} S_n\Rightarrow\int_0^{1}X_{\mathbf C}(s)^TX_{\mathbf C}(s)
\,ds.
\end{equation}
\end{theorem}

\begin{Remark*}
Theorem~\ref{theo4} is similar to Theorem 7 of Zhang and Chen \cite
{Zhang2007}.
Both characterize the asymptotic distribution of $S_n$. In
particular,
Zhang and Chen~\cite{Zhang2007} delineate the distribution of $\int_0^{1}X_{\mathbf C}(s)^TX_{\mathbf C}(s) \,ds$
as a $\chi^2$-type mixture.
All discussions associated with Theorem 7 of Zhang and Chen \cite
{Zhang2007} are valid
here, and therefore, we do not repeat them for the sake of space.

We consider conditional convergence for bootstrapped stochastic
processes. We focus on the bootstrapped process $\{G_{j}(s)^{(g)}\dvtx
s\in[0, 1]\}$ as the arguments for establishing the wild bootstrap
method for approximating the null distribution of $S_n$ and the
bootstrapped process $\{G_{j}(s)^{(g)}\dvtx  s\in[0, 1]\}$ are similar.
\end{Remark*}
%
\begin{theorem}\label{theo5}
If Assumptions~\ref{assumC1}--\ref{assumC9a} are true,
then $G_{j}(s)^{(g)}(s)$ converges weakly to $G_j(s)$ conditioning on
the data, where $G_j(s)$ is a centered Gaussian process indexed by
$s\in[0, 1]$.
\end{theorem}
\begin{Remark*}
Theorem~\ref{theo5} validates the bootstrapped process of
$G_{j}(s)^{(g)}$. An interesting observation is that the bias
correction for $\hat B_j(s)$ in constructing $G_{j}(s)^{(g)}$
is unnecessary. It leads to substantial computational saving.
\end{Remark*}

\section{Simulation studies}\label{sec5}

In this section, we present two simulation example to demonstrate the
performance
of the proposed procedures.

\begin{example}\label{example1} This example is designed to evaluate the
type I error rate and power of the proposed global test $S_n$
using Monte Carlo simulation. In this example, the data were
generated from a bivariate MVCM as follows:
%
\begin{equation}
\label{tract3eqsim1} y_{ij}(s_m)={\mathbf x}_i^TB_j(s_m)+
\eta_{ij}(s_m)+\varepsilon_{ij}(s_m)
\qquad\mbox{for }j=1, 2,
\end{equation}
where $s_m \sim U[0,1]$, $(\varepsilon_{i1}(s_m), \varepsilon_{i
2}(s_m))^T\sim N((0, 0)^T, S_\varepsilon(s_m)=\operatorname{diag}(\sigma_1^2,
\sigma_2^2))$ and ${\mathbf x}_i=(1, x_{i1}, x_{i2})$
for all $i=1,\ldots,n$ and
$m=1,\ldots, M$. Moreover,
$(x_{i1},\break x_{i2})^T\sim N((0, 0)^T, \operatorname{diag}(1-2^{-0.5},
1-2^{-0.5})+2^{-0.5}(1, 1)^{\otimes2})$ and
$\eta_{ij}(s)=\break\xi_{ij1}\psi_{j1}(s)+\xi_{ij2}\psi_{j2}(s)$, where
$\xi_{ijl} \sim N(0,\lambda_{jl})$
for $j=1,2$ and $l=1,2$. Furthermore, $s_m$, $(x_{i1}, x_{i2})$,
$\xi_{i11}$, $\xi_{i12}$, $\xi_{i21}$, $\xi_{i22}$, $\varepsilon_{i1}(s_m)$, and
$\varepsilon_{i2}(s_m)$ are independent random variables.
We set $(\lambda_{11}, \lambda_{12},\sigma_1^2,
\lambda_{21},\lambda_{22},\sigma_2^2)=(1.2, 0.6, 0.2, 1, 0.5, 0.1)$ and
the functional coefficients
and eigenfunctions as follows:
\begin{eqnarray*}
b_{11}(s)&=&s^2, b_{12}(s)=(1-s)^2,\qquad
b_{13}(s)=4s(1-s)-0.4;
\\
\psi_{11}(s)&=&\sqrt{2}\sin(2\pi s),\qquad \psi_{12}(s)=\sqrt{2}
\cos(2\pi s);
\\
b_{21}(s)&=&5(s-0.5)^2,\qquad b_{22}(s)=s^{0.5},\qquad
b_{23}(s)=4s(1-s)-0.4;
\\
\psi_{21}(s)&=&\sqrt{2}\cos(2\pi s),\qquad \psi_{22}(s)=\sqrt{2}
\sin(2\pi s).
\end{eqnarray*}
Then,
except for $( b_{13}(s)$, $ b_{23}(s))$ for all $s$,
we fixed all other parameters at the values specified above, whereas
we assumed $( b_{13}(s), b_{23}(s))=c(4s(1-s)-0.4$, $4s(1-s)-0.4)$,
where $c$ is a scalar specified below.

We want to test the hypotheses $H_0\dvtx  b_{13}(s)= b_{23}(s)=0$ for
all $s$ against $H_1\dvtx
b_{13}(s)\not=0$ or $b_{23}(s)\not=0$ for at least one
$s$. We set $c=0$ to assess the type~I error
rates for $S_n$, and set $c=0.1, 0.2, 0.3$ and $0.4$ to examine the
power of $S_n$. We set $M=50$, $n=200$ and $100$. For each
simulation, the significance levels were set at $\alpha= 0.05$ and
$0.01$, and 100 replications were used to estimate the rejection
rates.

\begin{figure}

\includegraphics{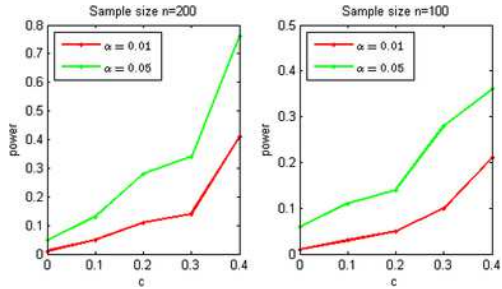}

\caption{Plot of power curves. Rejection rates of $S_n$ based on the
wild bootstrap method are calculated at five different values of $c$
(0, 0.1, 0.2, 0.3, and 0.4) for two sample sizes of $n$ (100 and
200) subjects at $5\%$ (green) and $1\%$ (red) significance levels.}
\label{tract3fig1}
\end{figure}

Figure~\ref{tract3fig1} depicts the power curves. It can be seen
from Figure~\ref{tract3fig1} that the rejection rates for $S_n$ based
on the wild bootstrap method are accurate for moderate sample sizes,
such as ($n=100$ or 200) at both significance levels ($\alpha=0.01$
or 0.05). As expected, the power increases with the sample size.
\end{example}

\begin{figure}[b]

\includegraphics{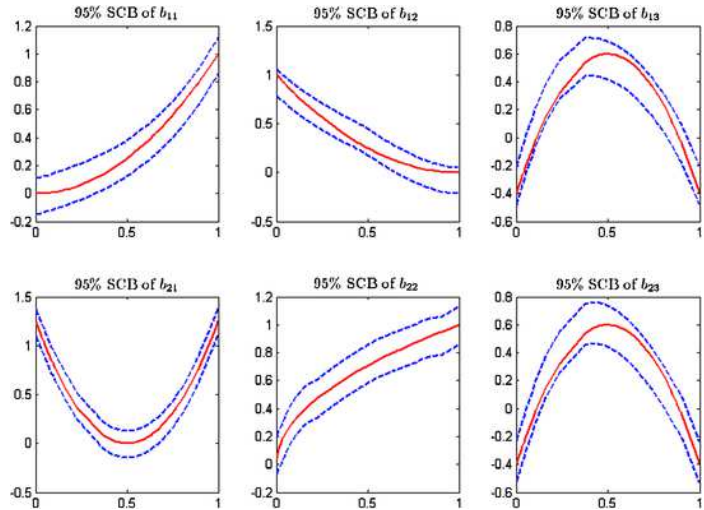}

\caption{Typical simultaneous confidence bands with $n=500$ and
$M=50$. The red solid curves are the true coefficient functions, and
the blue dashed curves are the confidence bands.}
\label{tract3fig2}
\end{figure}

\begin{example}\label{example2}
This example is used to evaluate the
coverage probabilities of SCB of the functional coefficients ${\mathbf
B}(s)$ based on the wild bootstrap method. The data were
generated from model (\ref{tract3eqsim1}) under the same
parameter values. We set $n=500$
and $M=25$, $50$, and $75$ and generated 200 datasets for each combination.
Based on the generated data, we calculated SCB for each component
%
\begin{table}
\tablewidth=250pt
\caption{Empirical coverage probabilities of $1-\alpha$ SCB for all
components of $B_1(\cdot)$ and $B_2(\cdot)$ based on $200$ simulated
data sets}\label{tabl1}
\begin{tabular*}{\tablewidth}{@{\extracolsep{\fill}}lcccccc@{}}
\hline
$\bolds{M}$ & $\bolds{b_{11}}$ & $ \bolds{b_{12}}$ & $\bolds{ b_{13}}$
& $\bolds{b_{21}}$ & $ \bolds{b_{22}}$ & $
\bolds{b_{23}}$\\
\hline
&\multicolumn{6}{c@{}}{$\alpha=0.05$}\\[4pt]
$25$ &0.915&0.930&0.945& 0.920&0.915&0.945\\
$50$ &0.925&0.940&0.945& 0.930&0.925&0.950 \\
$75$ &0.945&0.950&0.955& 0.945&0.945&0.955 \\
[6pt]
&\multicolumn{6}{c@{}}{$\alpha=0.01$}\\[4pt]
$25$ &0.985&0.965&0.985& 0.985&0.990&0.980\\
$50$ &0.995&0.980&0.985& 0.985&0.995&0.985 \\
$75$ &0.990&0.985&0.990& 0.995&0.990&0.990 \\
\hline
\end{tabular*}
\end{table}
of $B_1(s)$ and $B_2(s)$. Table~\ref{tabl1} summarizes the empirical coverage
probabilities based on $200$ simulations for $\alpha=0.01$ and
$\alpha=0.05$. The coverage probabilities improve with the number
of grid points $M$. When $M=75$, the differences between the
coverage probabilities and the claimed confidence levels are fairly
acceptable. The Monte Carlo errors are of size $\sqrt{0.95 \times
0.05/200}\approx0.015$ for $\alpha=0.05$. Figure~\ref{tract3fig2}
depicts typical simultaneous confidence bands, where $n=500$ and
$M=50$. Additional simulation results are given in the supplemental
article~\cite{ZLK2012}.
\end{example}



\section{Real data analysis}\label{sec6}


The data set consists of 128 healthy infants ($75$ males and $53$ females)
from the neonatal project on early brain
development.
The gestational ages of these infants
range from 262 to 433 days, and their mean gestational age
is 298 days with standard deviation 17.6 days.
The DTIs and T1-weighted images were acquired for each subject. For the
DTIs, the imaging parameters were as follows:
the six noncollinear directions at
the $b$-value of 1000 s/mm$^2$ with a
reference scan ($b=0$), the isotropic voxel $\mbox{resolution}=2$ mm, and the
in-plane field of $\mbox{view}=256$ mm in both
directions. A
total of five repetitions were acquired to improve the signal-to-noise
ratio of the DTIs.

The DTI data were processed by two key steps including a weighted least squares
estimation method~\cite{Basser1994b,Zhu2007b} to construct the
diffusion tensors and
a DTI atlas building pipeline~\cite{Goodlett2009,Zhu2010} to register
DTIs from multiple subjects to create
a study specific unbiased DTI atlas, to track fiber tracts in the atlas space
and to propagate them back into each subject's native space by using
registration information.
Subsequently, diffusion tensors (DTs) and their scalar diffusion
properties were calculated at each location along each individual fiber
tract by using
DTs in neighboring voxels close to the fiber tract.
Figure~\ref{tract3fig0}(a) displays the fiber bundle of the genu of
the corpus callosum (GCC), which is an area of white matter in the brain.
The GCC is the anterior end of the corpus callosum, and is bent
downward and backward in front of the septum pellucidum; diminishing
rapidly in thickness, it is prolonged backward under the name of the
rostrum, which is connected below with the lamina terminalis. It was
found that
neonatal microstructural development of GCC positively correlates with
age and callosal thickness.\looseness=1

The two aims of this analysis are to compare diffusion properties
including FA and MD along the GCC between the male and female groups
and to delineate the development of fiber diffusion properties across
time, which is addressed by including the gestational age at MRI
scanning as a covariate. FA and MD, respectively, measure the
inhomogeneous extent of local barriers to water diffusion and the
averaged magnitude of local water diffusion.
We fit model (\ref{Tract3eq1}) to the FA and MD values from all 128
subjects, in which
${\mathbf x}_i=(1, \mathrm{G}, \mathrm{Age})^T$, where $\mathrm{G}$ represents
gender. We then applied the estimation and inference procedures to
estimate ${\mathbf B}(s)$ and calculate $S_n$ for each hypothesis test.
We approximated the $p$-value of $S_n$ using the
wild bootstrap method with $G=1000$ replications. Finally, we
constructed the $95\%$ simultaneous confidence bands for the
functional coefficients of $B_j(s)$ for $j=1,2$.

\begin{figure}
\begin{tabular}{@{}c@{\hspace*{6pt}}c@{}}

\includegraphics{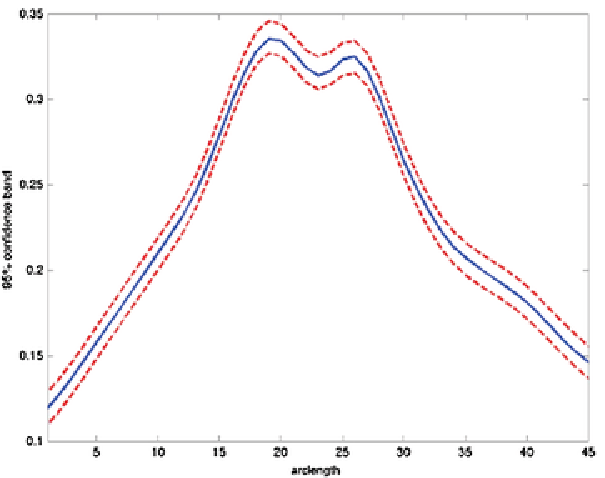}
 & \includegraphics{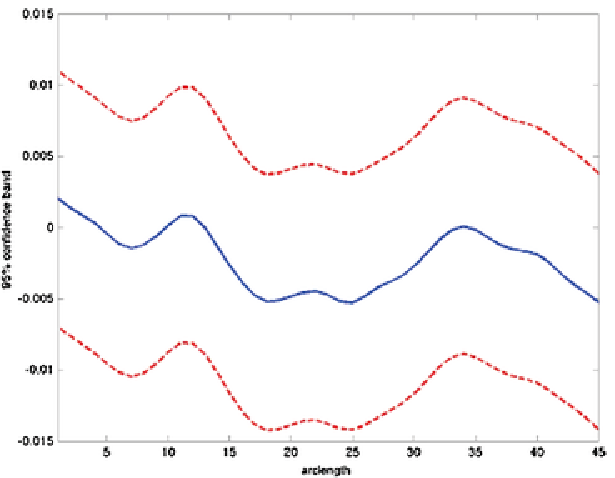}\\
(a) & (b)\\[6pt]

\includegraphics{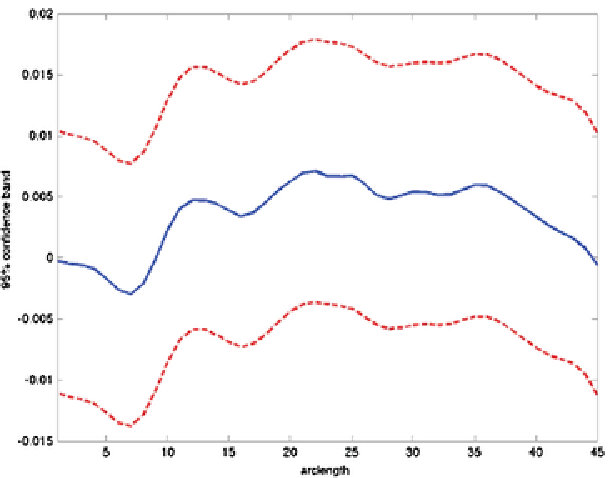}
 & \includegraphics{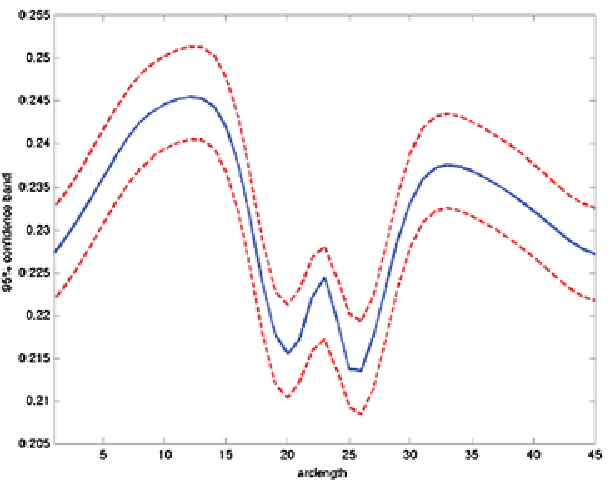}\\
(c) & (d)\\[6pt]

\includegraphics{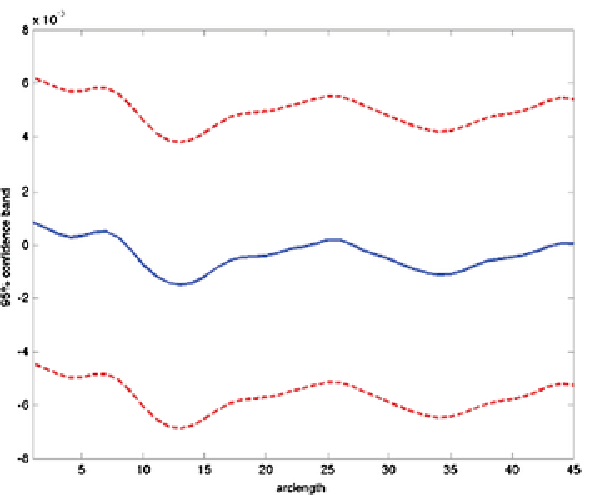}
 & \includegraphics{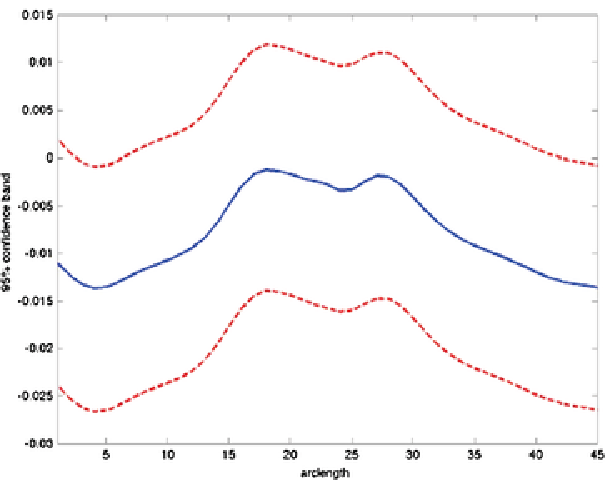} \\
(e) & (f)
\end{tabular}
\caption{Plot of estimated effects of intercept [\textup{(a)},
\textup{(d)}], gender [\textup{(b)}, \textup{(e)}], and age
[\textup{(c)}, \textup{(f)}] and their $95\%$ confidence bands. The
first three panels [\textup{(a)}, \textup{(b)}, \textup{(c)}] are for
FA and the last three panels [\textup{(d)},~\textup{(e)} and
\textup{(f)}] are for MD. The blue solid curves are the estimated
coefficient functions, and the red dashed curves are the confidence
bands.} \label{tract3fig6}
\end{figure}

Figure~\ref{tract3fig6} presents the estimated coefficient functions
corresponding to 1, G and Age associated with FA and MD (blue solid
lines in all panels of Figure~\ref{tract3fig6}).
The intercept functions [panels (a) and (d) in Figure~\ref{tract3fig6}] describe the
overall trend of FA and MD.
The gender coefficients for FA and MD
in Figure~\ref{tract3fig6}(b) and (e) are negative
at most of the grid points, which may indicate that compared with
female infants, male infants have relatively smaller magnitudes of
local water diffusivity along the genu of the corpus callosum.
The gestational
age coefficients for FA [panel (c) of Figure~\ref{tract3fig6}] are
positive at most grid points, indicating that FA measures increase
with age in both male and female infants, whereas those
corresponding to MD [panel~(f) of Figure~\ref{tract3fig6}] are
negative at most grid points.
This may indicate a negative correlation between the magnitudes of
local water diffusivity and gestational age along the genu of
the corpus callosum.

We statistically tested the effects of gender and gestational age on
FA and MD along the GCC tract. To test the
gender effect, we
computed the global test statistic $S_n=144.63$ and its associated
$p$-value ($p=0.078$), indicating a weakly significant gender
effect, which agrees with the findings in panels (b) and
(e) of Figure~\ref{tract3fig6}. A
moderately significant age effect was found with $S_n=929.69$
($p\mbox{-value}<0.001$). This agrees with the findings in panel (f) of
Figure~\ref{tract3fig6}, indicating that MD along the GCC tract
changes moderately with
gestational age. Furthermore, for FA and MD, we
constructed the $95\%$ simultaneous confidence bands of the
varying-coefficients for G$_i$ and age$_i$ (Figure~\ref{tract3fig6}).

\begin{figure}
\begin{tabular}{@{}cc@{}}

\includegraphics{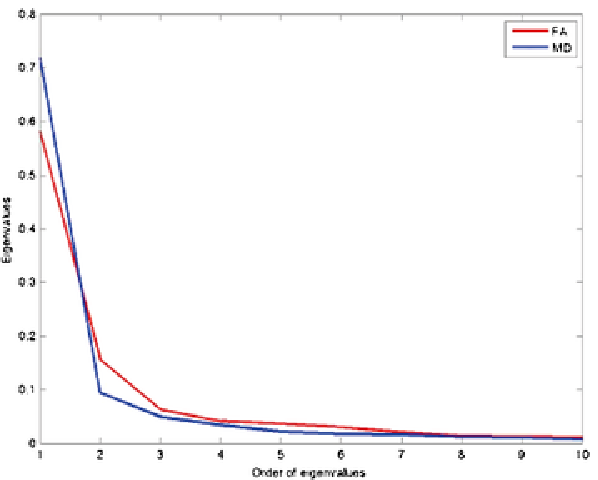}
 & \includegraphics{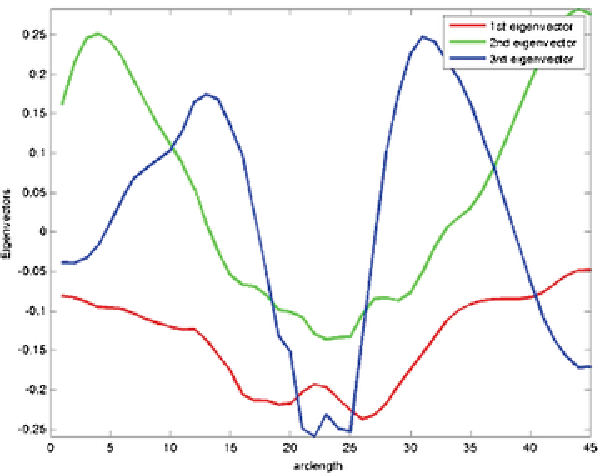}\\
(a) & (b)\\[6pt]
\multicolumn{2}{@{}c@{}}{
\includegraphics{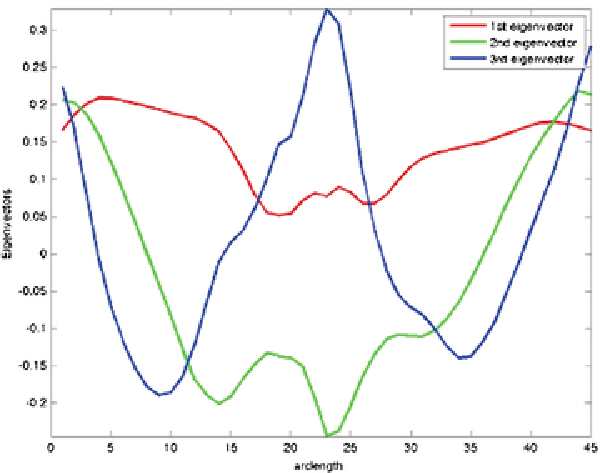}
}\\
\multicolumn{2}{@{}c@{}}{(c)}
\end{tabular}
\caption{Plot of the first $10$ eigenvalues \textup{(a)} and the first $3$
eigenfunctions for FA \textup{(b)} and MD \textup{(c)}.}
\label{tract3fig7}
\end{figure}

Figure~\ref{tract3fig7} presents the first $10$ eigenvalues and $3$
eigenfunctions of $\hat{\Sigma}_{\eta,{jj}}(s,t)$ for $j=1,2$. The
relative eigenvalues of $\hat\Sigma_{\eta, jj}$ defined as the ratios
of the eigenvalues of $\hat{\Sigma}_{\eta,{jj}}(s,t)$ over their sum
have similar distributional patterns [panel~(a) of Figure
\ref{tract3fig7}]. We observe that the first three eigenvalues account
for more than $90\%$ of the total and the others quickly vanish to
zero. The eigenfunctions of FA corresponding to the largest three
eigenvalues [Figure~\ref{tract3fig7}(b)] are different from those of MD
[Figure~\ref{tract3fig7}(c)].

In the supplement article~\cite{ZLK2012}, we further illustrate the
proposed methodology by an
empirical analysis of another real data set.




\begin{appendix}\label{app}
\section*{Appendix}

We introduce some notation. We define
%
\begin{eqnarray}\label{NewEq0}
T_{B, j} (h, s)&=&\sum_{i=1}^n
\sum_{m=1}^{M} K_{h}(s_m-s)
\bigl[{\mathbf x}_i\otimes {\mathbf z}_{h}(s_m-s)
\bigr] {\mathbf x}_i^TB_j(s_m),
\nonumber
\\
T_{\eta, j}(h, s)&=&\sum_{i=1}^n
\sum_{m=1}^{M} K_{h}(s_m-s)
\bigl[{\mathbf x}_i\otimes {\mathbf z}_{h}(s_m-s)
\bigr] \eta_{ij} (s_m),
\nonumber\\
T_{\varepsilon, j}(h, s)&=&\sum_{i=1}^n
\sum_{m=1}^{M} K_{h}(s_m-s)
\bigl[{\mathbf x}_i\otimes {\mathbf z}_{h}(s_m-s)
\bigr] \varepsilon_{ij}(s_m),
\\
r_u(K; s, h)&=&\frac{u_2(K; s, h)^2-u_1(K; s, h)u_3(K; s,
h)}{u_0(K; s, h)u_2(K; s, h)-u_1(K; s, h)^2},
\nonumber
\\
H_{h}(s_m-s)&=&K_h(s_m-s){\mathbf
z}_h(s_m-s),
\nonumber
\\
\Delta_{j}(s; {\ve\eta}_i, h_{1j}) &=&
M^{-1}\sum_{m=1}^{M}
H_{h_{1j}}(s_m-s) \eta_{ij}(s_m)\nonumber\\
&&{}- \int
_0^1 H_{h_{1j}}(u-s)
\eta_{ij}(u) \pi(u)\,du, \nonumber
\end{eqnarray}
where $u_r(K; s,
h)= \int_0^{1} h^{-r} (u-s)^r K_h(u-s)\,du$ for $r\geq0$.
Throughout the proofs, $C_k$'s stand for a generic constant, and it
may vary from line to line.

The proofs of Theorems~\ref{theo1}--\ref{theo5} rely on the following lemmas whose proofs are
given in the supplemental article~\cite{ZLK2012}.
%
\begin{lemma}\label{lemma1} Under Assumptions~\ref{assumC1},~\ref{assumC3}--\ref{assumC5}
and~\ref{assumC7},
we have that for
each $j$,
%
\begin{equation}
\label{Lemma1Eq0}\quad \sup_{s\in[0, 1]} n^{-1/2} h_{1j}\bigl|T_{\varepsilon, j}(h_{1j},
s)\bigr| =O_p\bigl(\sqrt{Mh_{1j}|{\log h_{1j}}|}\bigr)=o_p(Mh_{1j}).
\end{equation}
\end{lemma}

\begin{lemma}\label{lemma2}
Under Assumptions~\ref{assumC1},~\ref{assumC4},~\ref{assumC5} and~\ref{assumC7},
we have
that for any $r\geq0$ and $j$,
\begin{eqnarray*}
\sup_{s\in[0, 1]}\biggl\llvert \int K_{h_{1j}}(u-s)\frac{(u-s)^r}{h_{1j}^r}\,d
\bigl[\Pi_{M}(u)-\Pi(u)\bigr]\biggr\rrvert &=&O_p
\bigl((Mh_{1j})^{-1/2}\bigr),
\\
\sup_{s\in[0, 1]}\biggl\llvert \int K_{h_{1j}}(u-s)\frac{(u-s)^r}{h_{1j}^r}
\varepsilon_{ij}(u)\,d\Pi_{M}(u)\biggr\rrvert &=&O_p
\bigl(({Mh_{1j}})^{-1/2}\sqrt{|{\log h_{1j}}|}\bigr),
\end{eqnarray*}
where $\Pi_{M}(\cdot)$ is the sampling distribution function
based on ${\mathcal S}=\{s_1,\ldots, s_{M}\}$,
and $\Pi(\cdot)$ is the distribution function of $s_m$.
\end{lemma}


%
\begin{lemma}\label{lemma3}
Under Assumptions~\ref{assumC2}--\ref{assumC5}, we have
%
\begin{equation}
\label{Lemma3Eq1} \sup_{s\in[0, 1]}\Biggl| n^{-1/2}\sum
_{i=1}^n {\mathbf x}_i\otimes
\Delta_{j}(s;{\ve\eta}_i, h_{1j}) \Biggr|
=o_p(1).
\end{equation}
\end{lemma}


\begin{lemma}\label{lemma4}
If Assumptions~\ref{assumC1} and~\ref{assumC3}--\ref{assumC6} hold,
then we have
%
\begin{eqnarray}\label{Lemma4Eq1}
\mathrm{E}\bigl[\hat{B}_j(s)|{\mathcal S}\bigr]-{B}_j(s)
&=& 0.5 h_{1j}^{2}u_2(K) \ddot B_j(s)
\bigl[1+o_p(1)\bigr],
\nonumber\\[-8pt]\\[-8pt]
\operatorname{Var}\bigl[\hat B_j(s)|{\mathcal S}\bigr] &=& n^{-1}
\Sigma_{\eta, jj}(s, s)\Omega_X^{-1}
\bigl[1+o_p(1)\bigr],
\nonumber
\end{eqnarray}
where $e_n(s)=O_p((M h_{1j})^{-1/2})$
with $E[e_n(s)]=0$.
\end{lemma}
%
\begin{lemma}\label{lemma5}
If Assumptions~\ref{assumC1} and~\ref{assumC3}--\ref{assumC6} hold, then
for $s=0$ or $1$, we have
%
\begin{eqnarray}
\label{Lemma5Eq1} \mathrm{E}\bigl[\hat B_j(s)|{\mathcal S}
\bigr]-B_j(s) &=& 0.5 h_{1j}^{2}r_u(K;
s, h_{1j}) \ddot B_j(s) \bigl[1+o_p(1)\bigr],
\nonumber\\[-8pt]\\[-8pt]
\operatorname{Var}\bigl[\hat B_j(s)|{\mathcal S}\bigr]&=& n^{-1}
\Sigma_{\eta, jj}(s, s)\Omega_X^{-1}
\bigl[1+o_p(1)\bigr].
\nonumber
\end{eqnarray}
\end{lemma}
%
\begin{lemma}\label{lemma6}
Under Assumptions~\ref{assumC1}--\ref{assumC9a}, we have
\begin{eqnarray*}
\sup_{(s, t)}n^{-1} \Biggl|\sum
_{i=1}^n\overline\varepsilon_{ij}(s)
\eta_{ij}(t)\Biggr|&=&O_p\bigl(n^{-1/2}(\log
n)^{1/2}\bigr),
\\
\sup_{(s, t)}n^{-1} \Biggl|\sum_{i=1}^n
\overline\varepsilon_{ij}(s)\Delta\eta_{ij}(t)\Biggr|&=&O_p
\bigl(n^{-1/2}(\log n)^{1/2}\bigr),
\\
\sup_{s}n^{-1} \Biggl|\sum_{i=1}^n
\overline\varepsilon_{ij}(s) {\mathbf x}_{i}\Biggr|&=&O_p
\bigl(n^{-1/2}(\log n)^{1/2}\bigr),
\\
\sup_{s}n^{-1} \Biggl|\sum_{i=1}^n
\Delta\eta_{ij}(s) {\mathbf x}_{i}\Biggr|&=&O_p
\bigl(n^{-1/2}(\log n)^{1/2}\bigr).
\end{eqnarray*}
\end{lemma}
%
\begin{lemma}\label{lemma7}
Under Assumptions~\ref{assumC1}--\ref{assumC9a}, we have
\[
\sup_{(s, t)}n^{-1} \Biggl|\sum
_{i=1}^n\overline \varepsilon_{ij}(s)
\overline\varepsilon_{ij}(t)\Biggr|=O\bigl((Mh_{2j})^{-1} + (
\log n/n)^{1/2}\bigr)=o_p(1).
\]
\end{lemma}

We present only the key steps in the proof of Theorem~\ref{theo1} below. 

\begin{pf*}{Proof of Theorem~\ref{theo1}}
Define
\begin{eqnarray*}
{\mathbf U}_2(K; s, {\mathbf H})&=&\operatorname{diag}\bigl(r_u(K; s,
h_{11}),\ldots, r_u(K; s, h_{1J})\bigr),
\\
X_n(s)&=&\sqrt{n}\bigl\{\hat{\mathbf B}(s)-E\bigl[\hat{\mathbf B}(s)|{
\mathcal S}\bigr]\bigr\},
\\
X_{n, j}(s)&=&\sqrt{n}\bigl\{\hat
B_j(s)-E\bigl[\hat B_j(s)|{\mathcal S}\bigr]\bigr\}.
\end{eqnarray*}
According to the definition of
$\operatorname{vec}(\hat A_j(s))$, it is easy to see that
%
\begin{eqnarray}\qquad
\operatorname{vec}\bigl(\hat A_j(s)\bigr)&=& \Sigma(s,
h_{1j})^{-1} \bigl[T_{B, j}(h_{1j},
s)+T_{\varepsilon, j}(h_{1j}, s)+T_{\eta, j}(h_{1j}, s)
\bigr],
\\
X_{n, j}(s)&=& \sqrt{n} \bigl[{\mathbf I}_p\otimes(1, 0)
\bigr]\Sigma(s, h_{1j})^{-1}\bigl[T_{\varepsilon, j}(h_{1j},
s)+T_{\eta, j}(h_{1j}, s)\bigr].
\end{eqnarray}

The proof of Theorem~\ref{theo1}(i) consists of two parts:
\begin{itemize}
\item Part 1 shows that $\sqrt{n} \Sigma(s, h_{1j})^{-1}T_{\varepsilon,
j}(h_{1j}, s)=o_p(1)$ holds uniformly for all $s\in[0, 1]$ and $j=1,\ldots, J$.
\item
Part 2 shows that $\sqrt{n} \Sigma(s, h_{1j})^{-1}T_{\eta,
j}(h_{1j}, s)$ converges weakly to a Gaussian process $G(\cdot)$ with
mean zero and covariance matrix $\Sigma_{\eta, jj}(s, s')\Omega_X^{-1}$ for each $j$.
\end{itemize}

In part 1, we show that
%
\begin{equation}
\label{NewEq1} \sqrt{n} \bigl[{\mathbf I}_p\otimes(1, 0)\bigr]\Sigma(s,
h_{1j})^{-1}T_{\varepsilon,
j}(h_{1j}, s) =
o_p(1).
\end{equation}
It follows from Lemma~\ref{lemma1} that
\[
{n^{-1/2}} \sum_{i=1}^n {\mathbf
x}_i\otimes\Biggl\{M^{-1}\sum
_{m=1}^{M} K_{h_{1j}}(s_m-s) {\mathbf
z}_{h_{1j}}(s) \varepsilon_{i, j}(s_m)\Biggr\}
=o_p(1)
\]
hold uniformly for all $s\in[0, 1]$.
It follows from Lemma~\ref{lemma2} that
%
\begin{equation}\label{NewEq2}
(nM)^{-1}\Sigma(s, h_{1j}) = \Omega_X\otimes
\Omega_1(h_{1j}, s)+o_p(1)
\end{equation}
hold uniformly for all $s\in[0, 1]$. Based on these results, we can
finish the proof of
(\ref{NewEq1}).

In part 2, we show the weak convergence of $\sqrt{n}
[{\mathbf I}_p\otimes(1, 0)]\Sigma(s, h_{1j})^{-1}\*T_{\eta, j}(h_{1j},
s)$ for $j=1,\ldots, J$.
Part 2 consists of two steps. In Step~\ref{step1}, it follows from the standard
central limit theorem
that for each $s\in[0, 1]$,
%
\begin{equation}\label{proofeq7}\qquad
\sqrt{n} \bigl[{\mathbf I}_p\otimes(1, 0)\bigr]\Sigma(s,
h_{1j})^{-1}T_{\eta,
j}(h_{1j}, s)
\rightarrow^{L} N\bigl({\mathbf0}, \Sigma_{\eta, jj}(s, s)
\Omega_X^{-1}\bigr),
\end{equation}
where $\rightarrow^L$ denotes convergence in distribution.

Step~\ref{step2} shows the asymptotic tightness of $ \sqrt{n} [{\mathbf I}_p\otimes
(1, 0)]\Sigma(s, h_{1j})^{-1}\*T_{\eta, j}(h_{1j}, s)$. By using (\ref
{NewEq2}) and (\ref{NewEq0}), $\sqrt{n} \Sigma(s,
h_{1j})^{-1}T_{\eta, j}(h_{1j}, s)[1+o_p(1)] $
can be approximated by the sum of three terms (I), (II) and (III) as follows:
%
\begin{eqnarray}\label{NEWZHEQ5}\qquad
\mbox{(I)}&=& {n^{-1/2}} \sum_{i=1}^n
\Omega_X^{-1} {\mathbf x}_i\otimes
\Omega_1(h_{1j}, s)^{-1} \Delta_{j}(s;
{\ve\eta}_i, h_{1j}),
\nonumber\\
\mbox{(II)}&=&{n^{-1/2}} \sum_{i=1}^n
\Omega_X^{-1} {\mathbf x}_i\otimes
\Omega_1(h_{1j}, s)^{-1} \eta_{ij}(s)\nonumber\\
&&{}\times
\int_{\max(-sh_{1j}^{-1}, -1)}^{\min((1-s)h_{1j}^{-1}, 1)} K(u) (1, u)^T
\pi(s+h_{1j}u)\,du,
\nonumber\\[-8pt]\\[-8pt]
\mbox{(III)}&=&{n^{-1/2}} \sum_{i=1}^n
\Omega_X^{-1} {\mathbf x}_i\otimes
\Omega_1(h_{1j}, s)^{-1}
\nonumber
\\
&&\hspace*{40pt}{}\times \int_{\max(-sh_{1j}^{-1}, -1)}^{\min((1-s)h_{1j}^{-1}, 1)} K(u) \pmatrix{1
\cr
u} \bigl[
\eta_{ij}(s+h_{1j}u)-\eta_{ij}(s)\bigr]\nonumber\\
&&\hspace*{40pt}\qquad\quad\hspace*{51pt}{}\times
\pi(s+h_{1j}u)\,du.
\nonumber
\end{eqnarray}

We investigate the three terms on the
right-hand side of (\ref{NEWZHEQ5}) as follows.
It follows from Lemma~\ref{lemma3} that the first term on the
right-hand side of (\ref{NEWZHEQ5}) converges to zero uniformly.
We prove the asymptotic tightness of (II) as follows. Define
\begin{eqnarray*}
\hat X_{n, j}(s)&=& {n^{-1/2}} \sum_{i=1}^n
\Omega_X^{-1} {\mathbf x}_i\otimes(1, 0)
\Omega_1(h_{1j}, s)^{-1}\eta_{ij}(s)\\
&&\hspace*{0pt}{}\times
\int_{\max(-sh_{1j}^{-1}, -1)}^{\min((1-s)h_{1j}^{-1}, 1)} K(u) (1, u)^T
\pi(s+h_{1j}u)\,du.
\end{eqnarray*}
Thus,
we only need to prove the asymptotic tightness of $\hat X_{n, j}(s)$.
The asymptotic tightness of $\hat X_{n, j}(s)$ can be proved using the
empirical process techniques~\cite{VaartWellner1996}.
It follows that
\begin{eqnarray*}
&& (1, 0) \Omega_1(h_{1j}, s)^{-1} \int
_{\max(-sh_{1j}^{-1},
-1)}^{\min((1-s)h_{1j}^{-1}, 1)} K(u) (1, u)^T\pi
(s+h_{1j}u)\,du
\\
&&\qquad= \frac{u_2(K; s, h_{1j})u_0(K; s, h_{1j})-u_1(K; s, h_{1j})^2
+o(h_{1j})}{u_2(K; s, h_{1j})u_0(K; s, h_{1j}) -u_1(K; s, h_{1j})^2+o(h_{1j})} = 1+o(h_{1j}).
\end{eqnarray*}
Thus, $\hat X_{n, j}(s)$ can be simplified as
\[
\hat X_{n, j}(s)= \bigl[1+o(h_{1j}) \bigr] {n^{-1/2}}
\sum_{i=1}^n \eta_{ij}(s)
\Omega_X^{-1} {\mathbf x}_i.
\]
We consider a function class $
{\mathcal E}_{\eta}=\{ f(s; {\mathbf x}, \eta_{\cdot, j})=
\Omega_X^{-1}{\mathbf x} \eta_{\cdot,j}(s)\dvtx  s\in[0, 1]\}
$. Due to Assumption~\ref{assumC2}, ${\mathcal E}_{\eta}$ is a $P$-Donsker
class.

Finally, we consider the third term (III) on the right-hand side of
(\ref{NEWZHEQ5}).
It is easy to see that (III) can be written as
\begin{eqnarray*}
&& \Omega_X^{-1} \otimes\Omega_1(h_{1j},
s)^{-1} \\
&&\qquad{}\times\int_{\max
(-sh_{1j}^{-1}, -1)}^{\min((1-s)h_{1j}^{-1}, 1)} K(u)
\Biggl[{n^{-1/2}} \sum_{i=1}^n{
\mathbf x}_i\bigl\{\eta_{ij}(s+h_{1j}u)-
\eta_{ij}(s)\bigr\}\Biggr]\otimes \pmatrix{1
\cr
u} \\
&&\hspace*{83.2pt}\qquad{}\times\pi(s+h_{1j}u)\,du.
\end{eqnarray*}
Using the same argument of proving the second term (II), we can show
the asymptotic tightness of
${n^{-1/2}} \sum_{i=1}^n{\mathbf x}_i\eta_{ij}(s)$. Therefore, for any
$h_{1j}\rightarrow0$,
%
\begin{equation}
\label{NewEq4} \sup_{s\in[0, 1], |u|\leq1} \Biggl|n^{-1/2} \sum
_{i=1}^n{\mathbf x}_i\bigl\{
\eta_{ij}(s+h_{1j}u)-\eta_{ij}(s)\bigr
\}\Biggr|=o_p(1).
\end{equation}
It follows from Assumptions~\ref{assumC5} and~\ref{assumC7} and
(\ref{NewEq4}) that (III) converges to zero uniformly. Therefore, we
can finish the proof of Theorem~\ref{theo1}(i). Since
Theorem~\ref{theo1}(ii) is a direct consequence of
Theorem~\ref{theo1}(i) and Lemma~\ref{lemma4}, we finish the proof of
Theorem~\ref{theo1}.
\end{pf*}

\begin{pf*}{Proof of Theorem~\ref{theo2}} Proofs of parts (a)--(d) are
completed by some straightforward calculations. Detailed
derivation is given in the supplemental document. Here we prove part
(e) only.
Let $\tilde K_{M, h}(s)=\tilde K_{M}(s/h)/h$, where
$\tilde K_{M}(s)$ is the empirical equivalent kernels for the
first-order local polynomial kernel~\cite{Fan1996}.
Thus, we have
%
\begin{eqnarray}\label{Th2Eq1}\quad
\hat\eta_{ij}(s)-\eta_{ij}(s) &=& \sum
_{m=1}^{M} \tilde K_{M,h_{2j}}(s_{m}-s){
\mathbf x}_i^T\bigl[B_j(s_m)- \hat
B_j(s_m)\bigr]
\nonumber\\[-9pt]\\[-9pt]
&&{} + \sum_{m=1}^{M} \tilde
K_{M,h_{2j}}(s_{m}-s)\bigl[\eta_{ij}(s_m)+
\varepsilon_{ij}(s_m)-\eta_{ij}(s)\bigr].\nonumber
\end{eqnarray}
We define
\begin{eqnarray*}
\overline\varepsilon_{ij}(s)&=& \sum
_{m=1}^{M} \tilde K_{M,h_{2j}}(s_{m}-s)
\varepsilon_{ij}(s_m),
\\[-1pt]
\Delta\eta_{ij}(s)&=& \sum_{m=1}^{M}
\tilde K_{M,h_{2j}}(s_{m}-s)\bigl[\eta_{ij}(s_m)-
\eta_{ij}(s)\bigr],
\\[-1pt]
\Delta B_{j}(s)&=&\sum_{m=1}^{M}
\tilde K_{M,h_{2j}}(s_{m}-s) \bigl[B_j(s_m)-
\hat B_{j}(s_m)\bigr],
\\[-1pt]
\Delta_{ij}(s)&=&\overline\varepsilon_{ij}(s)+\Delta
\eta_{ij}(s)+{\mathbf x}_i^T \Delta
B_{j}(s).
\end{eqnarray*}
It follows from (\ref{Th2Eq1}) that
%
\begin{equation}
\label{Th2Eq5} \hat\eta_{ij}(s)-\eta_{ij}(s)=
\Delta_{ij}(s)=\overline\varepsilon_{ij}(s)+\Delta
\eta_{ij}(s)+{\mathbf x}_i^T \Delta
B_{j}(s).
\end{equation}
It follows from Lemma~\ref{lemma2} and a Taylor expansion that
\[
\sup_{s\in[0, 1]} \bigl|\overline \varepsilon_{ij}(s)\bigr|=O_p
\biggl(\sqrt{\frac{|{\log(h_{2j})}|}{Mh_{2j}}}\biggr)
\]
and
\[
\sup_{s\in[0, 1]} \bigl|\Delta
\eta_{ij}(s)\bigr|=O_p(1)\sup_{s\in[0,
1]} \bigl|\ddot
\eta_{ij}(s)\bigr| h_{1j}^{(2) 2}.\vadjust{\goodbreak}
\]
Since $\sqrt{n} \{\hat B_{j}(\cdot)-B_j(\cdot)-0.5u_2(K)^2 h_{1j}^2\ddot B_j(\cdot)[1+o_p(1)]\}$ weakly converges to a Gaussian
process in $\ell^{\infty}([0, 1])$ as
$n\rightarrow\infty$, $\sqrt{n} \{\hat B_{j}(\cdot)-B_j(\cdot
)-0.5u_2(K)^2 h_{1j}^2\*\ddot B_j(\cdot)[1+o_p(1)]\}$ is asymptotically
tight. Thus, we have
\begin{eqnarray*}
\Delta B_{ij}(s)&=&-\sum_{m=1}^{M}
\tilde K_{M,h_{2j}}(s_{j}-s) 0.5u_2(K)^2
h_{1j}^2\ddot B_j(s_m)
\bigl[1+o_p(1)\bigr]
\\
&&{} +\sum_{m=1}^{M} \tilde
K_{M,h_{2j}}(s_{j}-s) \bigl\{0.5u_2(K)^2
h_{1j}^2\ddot B_j(s_m)
\bigl[1+o_p(1)\bigr]\\
&&\hspace*{160.2pt}{}+B_j(s_m)-\hat
B_{j}(s_m)\bigr\},
\end{eqnarray*}
\begin{eqnarray*}\\[-15pt]
\sup_{s\in[0, 1]} \bigl\| \Delta B_{ j}(s)\bigr\|&=&O_p
\bigl(n^{-1/2}\bigr)+O_p\bigl( h_{1j}^2
\bigr).
\end{eqnarray*}
Combining these results, we have
\[
\sup_{s\in[0, 1]} \bigl| \hat\eta_{ij}(s)-\eta_{ij}(s)\bigr|
=O_p\bigl({\bigl|\log (h_{2j})\bigr|}^{1/2}({Mh_{2j}})^{-1/2}+h_{1j}^{(2) 2}+
h_{1j}^2+n^{-1/2}\bigr).
\]
This completes the proof of part (e).
\end{pf*}

\begin{pf*}{Proof of Theorem~\ref{theo3}} Recall that
$ \hat\eta_{ij}(s)= \eta_{ij}(s)+\Delta_{ i, j}(s)$, we have
%
\begin{eqnarray}\label{Th3Eq1}\qquad
n^{-1}\sum_{i=1}^n
\hat\eta_{ij}(s) \hat\eta_{ij}(t)&=&n^{-1} \sum
_{i=1}^n\Delta_{ij}(s)
\Delta_{ij}(t)+ n^{-1} \sum_{i=1}^n
\eta_{ij}(s) \Delta_{ij}(t)
\nonumber\\[-8pt]\\[-8pt]
&&{}+n^{-1}\sum_{i=1}^n
\Delta_{ij}(s) \eta_{ij}(t)+n^{-1}\sum
_{i=1}^n\eta_{ij}(s) \eta_{ij}(t).\nonumber
\end{eqnarray}
This proof consists of two steps. The first step is to
show that the first three terms on the right-hand side of (\ref
{Th3Eq1}) converge to zero uniformly for all $(s, t)\in[0, 1]^2$ in
probability.
The second step is to show the uniform convergence of $n^{-1}\sum_{i=1}^n\eta_{ij}(s) \eta_{ij}(t)$ to $\Sigma_\eta(s, t)$ over
$(s, t)\in[0, 1]^2$ in probability.

We first show that
%
\begin{equation}
\label{Th3PEq1} \qquad\sup_{(s, t)}n^{-1} \Biggl|\sum
_{i=1}^n\Delta_{ij}(s)
\eta_{ij}(t)\Biggr|=O_p\bigl(n^{-1/2}+
h_{1j}^2+h_{2j}^{2}+(\log
n/n)^{1/2}\bigr).
\end{equation}
Since
%
\begin{eqnarray}
\label{Th3PEq2}
&&
\sum_{i=1}^n
\Delta_{ij}(s) \eta_{ij}(t) \nonumber\\
&&\qquad\leq n^{-1} \Biggl\{ \Biggl|
\sum_{i=1}^n\overline\varepsilon_{ij}(s)
\eta_{ij}(t)\Biggr|+ \Biggl|\sum_{i=1}^n
\Delta\eta_{ij}(s)\eta_{ij}(t)\Biggr|\\
&&\qquad\quad\hspace*{85pt}{}+ \Biggl|\sum
_{i=1}^n{\mathbf x}_i^T \Delta
B_{j}(s) \eta_{ij}(t)\Biggr|\Biggr\},\nonumber
\end{eqnarray}
it is sufficient to focus on the three terms on the right-hand side
of (\ref{Th3PEq2}). Since
\[
\bigl|{\mathbf x}_i^T \Delta B_{j}(s)
\eta_{ij}(t)\bigr| \leq\|{\mathbf x}_i\|_2
\sup_{s\in[0, 1]} \bigl\| \Delta B_{ k}(s)\bigr\|_2
\sup_{t\in[0, 1]}\bigl|\eta_{ij}(t)\bigr|,
\]
we have
\begin{eqnarray*}
n^{-1} \Biggl|\sum_{i=1}^n{\mathbf
x}_i^T \Delta B_{j}(s) \eta_{ij}(t)\Biggr|
&\leq& \sup_{s\in[0, 1]} \bigl\| \Delta B_{ k}(s)\bigr\|_2
n^{-1} \sum_{i=1}^n \|{\mathbf
x}_i\|_2 \bigl|\eta_{ij}(t)\bigr|\\
&=&O_p
\bigl(n^{-1/2}+ h_{1j}^2\bigr).
\end{eqnarray*}
Similarly, we have
\[
n^{-1} \Biggl|\sum_{i=1}^n\Delta
\eta_{ij}(s)\eta_{ij}(t)\Biggr|\leq n^{-1} \sum
_{i=1}^n \sup_{s, t\in[0, 1]}\bigl|\Delta
\eta_{ij}(s)\eta_{ij}(t)\bigr|=O_p
\bigl(h_{1j}^{(2) 2}\bigr)=o_p(1).
\]
It follows from Lemma~\ref{lemma6} that $\sup_{(s, t)}n^{-1} \{
|\sum_{i=1}^n\overline\varepsilon_{ij}(s)\eta_{ij}(t)|=O((\log
n/n)^{1/2})$. Similarly, we can show that $ \sup_{(s, t)}n^{-1}
|\sum_{i=1}^n\Delta_{ij}(t) \eta_{ij}(s)|=O_p(n^{-1/2}+
h_{1j}^2+h_{2j}^{2}+(\log n/n)^{1/2})$.\vspace*{1pt}

We can show that
%
\begin{equation}
\label{Th3PEq3} \sup_{(s,t)} \Biggl|n^{-1}\sum
_{i=1}^n\bigl[ \eta_{ij}(s)
\eta_{ij}(t)-\Sigma_{\eta, jj}(s, t)\bigr]\Biggr|=O_p
\bigl(n^{-1/2}\bigr).
\end{equation}
Note that
\begin{eqnarray*}
&& \bigl| \eta_{ij}(s_1) \eta_{ij}(t_1)-
\eta_{ij}(s_2) \eta_{ij}(t_2)\bigr|
\\
&&\qquad \leq2\bigl(|s_1-s_2|+|t_1-t_2|\bigr)
\sup_{s\in[0, 1]}\bigl|\dot \eta_{ij}(s)\bigr|\sup_{s\in[0, 1]}\bigl|
\eta_{ij}(s)\bigr|
\end{eqnarray*}
holds for any $(s_1, t_1)$ and $(s_2, t_2)$,
the functional class $\{ \eta_{j}(u) \eta_{j}(v)\dvtx  (u, v)\in[0, 1]^2\}
$ is a Vapnik and Cervonenkis (VC) class \cite
{VaartWellner1996,Kosorok2008}. Thus, it yields that (\ref{Th3PEq3})
is true.

Finally, we can show that
%
\begin{eqnarray}
\label{Th3PEq4}
&&
\sup_{(s, t)}n^{-1} \Biggl|\sum
_{i=1}^n\Delta_{ij}(s)
\Delta_{ij}(t)\Biggr|\nonumber\\[-8pt]\\[-8pt]
&&\qquad=O_p\bigl((Mh_{2j})^{-1}
+ (\log n/n)^{1/2}+ h^4_{j}+h_{1j}^{(2)4}
\bigr).\nonumber
\end{eqnarray}
With some calculations, for a positive constant $C_1$, we have
\begin{eqnarray*}
&&\Biggl|\sum_{i=1}^n\Delta_{ij}(s)
\Delta_{ij}(t)\Biggr|\\
&&\qquad\leq C_1 \sup_{(s,t)} \Biggl[ \Biggl|\sum
_{i=1}^n\overline\varepsilon_{ij}(s)
\overline\varepsilon_{ij}(t)\Biggr|+ \Biggl|\sum_{i=1}^n
\overline\varepsilon_{ij}(s) \Delta\eta_{ij}(t)\Biggr|
\\
&&\hspace*{34pt}\qquad\quad{}+ \Biggl|\sum_{i=1}^n\Delta
\eta_{ij}(t){\mathbf x}_i^T\Delta
B_j(s)\Biggr| +\Biggl|\sum_{i=1}^n
\overline\varepsilon_{ij}(s){\mathbf x}_i^T\Delta
B_j(t)\Biggr|
\\
&&\hspace*{34pt}\qquad\quad{}+ \Biggl|\sum_{i=1}^n\Delta
\eta_{ij}(s) \Delta\eta_{ij}(t)\Biggr|+\Biggl|\sum
_{i=1}^n{\mathbf x}_i^T\Delta
B_j(s)\Delta B_j(t){\mathbf x}_i\Biggr|\Biggr].
\end{eqnarray*}

It follows from Lemma~\ref{lemma7} that
\begin{eqnarray*}
\sup_{(s, t)}n^{-1} \Biggl|\sum_{i=1}^n
\overline\varepsilon_{ij}(s) \overline\varepsilon_{ij}(t)\Biggr|&=&O_p
\bigl((Mh_{2j})^{-1} + (\log n/n)^{1/2}\bigr),\\[-15pt]
\end{eqnarray*}
\begin{eqnarray*}
&& \sup_{(s, t)}n^{-1}\Biggl[\Biggl|\sum
_{i=1}^n\overline\varepsilon_{ij}(s) \Delta
\eta_{ij}(t)\Biggr|+ \Biggl|\sum_{i=1}^n
\Delta\eta_{ij}(t){\mathbf x}_i^T\Delta
B_j(s)\Biggr| +\Biggl|\sum_{i=1}^n
\overline\varepsilon_{ij}(s){\mathbf x}_i^T\Delta
B_j(t)\Biggr|\Biggr]
\\[2pt]
&&\qquad= O_p\bigl( (\log n/n)^{1/2}\bigr).
\end{eqnarray*}
Since $ \sup_{s\in[0, 1]}
|\Delta\eta_{ij}(s)|=C_2 \sup_{s\in[0, 1]} |\ddot\eta_{ij}(s)|h_{2j}^{2}, $
we have
\[
\sup_{(s, t)}n^{-1}
\Biggl|\sum_{i=1}^n\Delta\eta_{ij}(s) \Delta\eta_{ij}(t)\Biggr|=O\bigl(
h_{1j}^{(2)4}\bigr).
\]
Furthermore,
since $\sup_{s\in[0, 1]}\|\Delta{\mathbf B}(s)\|=O_p(n^{-1/2}+
h^2_{j})$, we have
\[
n^{-1}\Biggl|\sum_{i=1}^n{\mathbf
x}_i^T\Delta B_j(s)\Delta
B_j(t){\mathbf x}_i\Biggr|=O_p\bigl(n^{-1}+
h^4_{j}\bigr).
\]
Note that
the arguments for (\ref{Th3PEq1})--(\ref{Th3PEq4}) hold for $\hat
\Sigma_{\eta, jj'}(\cdot, \cdot)$ for any $j\not=j'$. Thus,
combining (\ref{Th3PEq1})--(\ref{Th3PEq4}) leads to Theorem~\ref{theo3}(i).

To prove Theorem~\ref{theo3}(ii), we follow the same arguments in Lemma 6 of
Li and Hsing~\cite{LiHsing2010}. For completion, we highlight several
key steps
below. We define
%
\begin{equation}
\label{Cor1PEq1} (\Delta\psi_{j, j}) (s)=\int_0^{1}
\bigl[\hat\Sigma_{\eta, jj}(s, t)-\Sigma_{\eta, jj}(s, t)\bigr]
\psi_{j,
j}(t)\,dt.
\end{equation}
Following Hall and Hosseini-Nasab~\cite{HallHosseini2006} and the
Cauchy--Schwarz inequality, we have
\begin{eqnarray*}
&&\biggl\{\int_0^{1} \bigl[\hat
\psi_{j, j}(s)-\psi_{j, j}(s)\bigr]^2\,ds\biggr
\}^{1/2}
\\
&&\qquad \leq C_2\biggl\{ \biggl[\int_0^{1}
(\Delta\psi_{j, j}) (s)^2\,ds\biggr]^{1/2} + \int
_0^{1}\int_0^{1}
\bigl[\hat\Sigma_{\eta, jj}(s, t)-\Sigma_{\eta, jj}(s, t)
\bigr]^2 \,ds \,dt \biggr\}
\\
&&\qquad\leq C_2\biggl\{\int_0^{1}\int
_0^{1}\bigl[\hat\Sigma_{\eta, jj}(s, t)-
\Sigma_{\eta, jj}(s, t)\bigr]^2 \,ds \,dt \biggr\}^{1/2}
\biggl\{\int_0^{1}\bigl[\psi_{j,
j}(t)
\bigr]^2\,dt\biggr\}^{1/2}
\\
&&\qquad\quad{} +\int_0^{1}\int_0^{1}
\bigl[\hat\Sigma_{\eta, jj}(s, t)-\Sigma_{\eta, jj}(s, t)
\bigr]^2 \,ds \,dt
\\
&&\qquad \leq C_3 \sup_{(s, t)\in[0, 1]^2}\bigl|\hat\Sigma_{\eta, jj}(s, t)-
\Sigma_{\eta, jj}(s, t)\bigr|,
\end{eqnarray*}
which yields Theorem~\ref{theo3}(ii)(a).

Using (4.9) in Hall, M{\"u}ller and Wang~\cite{MR2278365}, we have
\begin{eqnarray*}
&&|\hat\lambda_{j, j}-\lambda_{j, j}|
\\
&&\qquad\leq |\int_0^{1}\int_0^{1}
[\hat\Sigma_{\eta, jj}-\Sigma_{\eta, jj}](s, t)\psi_{j, j}(s)
\psi_{j, j}(t)\,ds\,dt\\
&&\qquad\quad{}+O\biggl(\int_0^{1} (
\Delta\psi_{j, j}) (s)^2\,ds\biggr)
\\
&&\qquad\leq C_4 \sup_{(s, t)\in[0, 1]^2}\bigl|\hat\Sigma_{\eta, jj}(s, t)-
\Sigma_{\eta, jj}(s, t)\bigr|,
\end{eqnarray*}
which yields Theorem~\ref{theo3}(ii)(b). This completes the proof.
\end{pf*}

\begin{pf*}{Proof of Theorem~\ref{theo5}}
The proof of Theorem~\ref{theo5} is given in the supplement
arctile~\cite{ZLK2012}.
\end{pf*}
\end{appendix}

\section*{Acknowledgments}

The authors are grateful to the Editor Peter B$\ddot{\mbox{u}}$hlmann,
the Associate Editor, and three anonymous referees for valuable
suggestions, which have greatly helped to improve our presentation.

\begin{supplement}
\stitle{Supplement to ``Multivariate varying coefficient model for
functional responses''}
\slink[doi]{10.1214/12-AOS1045SUPP} 
\sdatatype{.pdf}
\sfilename{aos1045\_supp.pdf}
\sdescription{This supplemental material includes the proofs of all
theorems and lemmas.}
\end{supplement}

%

\printaddresses

\end{document}